\documentclass[hidelinks,onefignum,onetabnum]{siamart220329}
       
\usepackage[utf8]{inputenc}    
\usepackage[T1]{fontenc}    
\usepackage[english]{babel}
\usepackage{enumerate}
\usepackage{graphicx} 
\usepackage{subcaption}								
\usepackage{booktabs}
\usepackage{mathtools}
\usepackage{amsmath,bm}
\usepackage{amssymb} 
\usepackage{amsfonts}
\usepackage{fancyhdr}
\usepackage{fourier} 
\usepackage{pgfplots}  
\usepackage{cancel}  
 
\usepackage{epstopdf}  
\ifpdf
\DeclareGraphicsExtensions{.eps,.pdf,.png,.jpg}
\else
\DeclareGraphicsExtensions{.eps}
\fi

\usepackage[square,sort,comma,numbers]{natbib}                               
\usepackage{tikz-cd}
 
\numberwithin{equation}{section}					
\numberwithin{figure}{section}						
\numberwithin{table}{section}						

\usepackage{listings}								
\usepackage[top=3cm, bottom=3cm,
			inner=3cm, outer=3cm]{geometry}			
\usepackage{eso-pic}								

\usepackage{float} 									
\usepackage{parskip}								
\usepackage{tcolorbox}

\newcommand{\RR}{\mathbb{R}} 
   
\newcommand{\NN}{\mathbb{N}}

\newcommand{\RT}{\mathbf{RT}}
\newcommand{\BDM}{\mathbf{BDM}}
 
\newcommand{\FF}{\mathnormal{F}}
\newcommand{\GG}{\mathnormal{G}} 

\newcommand{\permi}{\boldsymbol{\eta}}
\DeclareMathOperator{\dive}{div}

\newcommand{\permig}{\permi_{\Gamma}}
\newcommand{\jump}[1]{\llbracket #1 \rrbracket}  
\newcommand{\bfV}{\mathbf{V}} 
\newcommand{\bfQ}{\mathbf{Q}} 

\newcommand{\bff}{\mathbf{f}}    
\newcommand{\bfK}{\mathbf{K}} 
\newcommand{\bfH}{\mathbf{H}}   
\newcommand{\bfE}{\mathbf{E}} 
\newcommand{\bfL}{\mathbf{L}}   
\newcommand{\bfx}{\mathbf{x}} 
\newcommand{\boundp}{\partial\Omega_{p}} 
\newcommand{\bfu}{\pmb{u}} 
\newcommand{\bfv}{\pmb{v}}  
\newcommand{\bfh}{\pmb{h}}  
\newcommand{\bfw}{\pmb{w}}
\newcommand{\bfn}{\pmb{n}}
\newcommand{\boundu}{\partial\Omega_{\pmb{u}}}
\newcommand{\mesh}{\mathcal{T}_h}
\newcommand{\mcT}{\mathcal{T}}
\newcommand{\mcM}{\mathcal{M}}

\newcommand{\mcF}{\mathcal{F}}

\newcommand{\interf}{\Gamma}

\newcommand{\pB}{p_B}              
\newcommand{\uB}{u_B}

\newcommand{\kp}{\hat{k}}
      
\newcommand{\mcTh}{\mcT_{h,i}} 
\newcommand{\mcTT}{\mcT_{h,i}(T_L)} 
\newcommand{\mcTL}{\mcT_{L,i}} 
\newcommand{\mcFM}{\mcF_{h,i}(M_L)} 
\newcommand{\mcGh}{\mathcal{G}_{h}}

\newcommand{\tn}{|\mspace{-1mu}|\mspace{-1mu}|}
 
\newcommand{\restreinta}[1]{\mathclose{}|\mathopen{}_{#1}} 
 
\newcommand{\Vki}{\bfV_k{(\mcT_{h,i})}}

\newcommand{\Vk}{\bfV_{k,h}}            
\newcommand{\Qki}{Q_{\kp}{(\mcT_{h,i})}}       
\newcommand{\Qk}{Q_{\kp,h}}

\newcommand{\Af}{\mathcal{A}}          
\newcommand{\Bf}{\mathcal{B}}   
\newcommand{\Cf}{\mathbf{C}} 
\newcommand{\intfneig}{U_{\varepsilon}(\Gamma)}   
              
\newcommand{\Om}{\Omega_{1}\cup  \Omega_{2}}    
\newcommand{\Omeps}{\Omega_{i,\varepsilon}}
  
\newcommand{\Omdhi}{\Omega_{\mcT_{h,i}}}
\newcommand{\Omdh}{\Omega_{\mcT_{h,1}} \cup  \Omega_{\mcT_{h,2}}}   
\newcommand{\interpv}{\pi_h^k}              
\newcommand{\interpvi}{\pi_{h,i}^k}  
 
\newcommand{\interpp}{\hat{\Pi}_{h}^{\hat{k}}}
  
\newcommand{\projpi}{\hat{\Pi}_{h,i}^{\hat{k}}}  
\newcommand{\projpii}[1]{\hat{\Pi}_{h,#1}^{\hat{k}}}   
\newcommand{\varepsilonz}{\varepsilon_0}                
                         
\newcommand{\intorp}{k+1}                   
\newcommand{\interr}{\pmb{\rho}} 
\newcommand{\bfrho}{\pmb{\rho}} 
\newcommand{\Hdivr}{\bfH^{r,\dive}}
\newcommand{\Hdivinto}{\bfH^{k+1,\dive}}     

\newsiamremark{assumption}{Assumption}           
\newsiamremark{remark}{Remark}         
\newsiamthm{thm}{Theorem}    
\newsiamthm{prop}{Proposition}        
  
\headers{A divergence preserving cut finite element method for Darcy flow}{T. Frachon, P. Hansbo, E. Nilsson, and S. Zahedi}

\title{A divergence preserving cut finite element method for Darcy flow \thanks{Submitted to the editors DATE.
		\funding{This research was supported by the Swedish Research Council Grant No. 2018-04192, 2022-04808, and the Wallenberg Academy Fellowship KAW 2019.0190.}}}

\author{Thomas Frachon\thanks{Department of Mathematics, KTH Royal Institute of Technology,  SE-100 44 Stockholm, Sweden 
		(\email{frachon@kth.se}, \email{erikni6@kth.se}, \email{sara.zahedi@math.kth.se}).}
	\and Peter Hansbo\thanks{Department of Mechanical Engineering, J\"onk\"oping University, SE-551~11 J\"onk\"oping, Sweden
		(\email{Peter.Hansbo@ju.se}).}
	\and Erik Nilsson \footnotemark[2]
		\and Sara Zahedi \footnotemark[2]
	}

\usepackage{amsopn}
\ifpdf
\hypersetup{
	pdftitle={A divergence preserving cut finite element method for Darcy flow},
	pdfauthor={T. Frachon, P. Hansbo, E. Nilsson, and S. Zahedi}
}
\fi

\begin{document}

\maketitle
 
\begin{abstract} 
We study cut finite element discretizations of a Darcy interface problem based on the mixed finite element pairs $\RT_{k} \times Q_k$, $k\geq 0$. Here $Q_k$ is the space of discontinuous polynomial functions of degree less or equal to $k$ and $\RT$ is the Raviart-Thomas space. We show that the standard ghost penalty stabilization, often added in the weak forms of cut finite element methods for stability and control of the condition number of the resulting linear system matrix, destroys the divergence-free property of the considered element pairs. Therefore, we propose new stabilization terms for the pressure and show that we recover the optimal approximation of the divergence without losing control of the condition number of the linear system matrix. We prove that with the new stabilization term the proposed cut finite element discretization results in pointwise divergence-free approximations of solenoidal velocity fields. We derive a priori error estimates for the proposed unfitted finite element discretization based on $\RT_k\times Q_k$, $k\geq 0$. In addition, by  decomposing the computational mesh into macro-elements and applying ghost penalty terms only on interior edges of macro-elements, stabilization is applied very restrictively and active only where needed.  Numerical experiments with element pairs $\RT_0\times Q_0$, $\RT_1\times Q_1$, and $\BDM_1\times Q_0$ (where $\BDM$ is the Brezzi-Douglas-Marini space) indicate that with the new method we have 1) optimal rates of convergence of the approximate velocity and pressure;  2) well-posed linear systems where the condition number of the system matrix scales as it does for fitted finite element discretizations; 3) optimal rates of convergence of the approximate divergence with pointwise divergence-free approximations of solenoidal velocity fields.  All three properties hold independently of how the interface is positioned relative to the computational mesh.
\end{abstract}      
                                   
                            
                        
\section{Introduction}\label{sec: intro}
We consider the Darcy interface problem \cite{DanSco12}, where a fractured porous medium is occupied by a fluid wherein the fracture (a thin region with different porous structure \cite{DanSco12}) is modelled by the interface. Several strategies of discretizing such problems exist, see \cite{FLe1BerBoo18} and references therein. Our aim is to develop unfitted finite element discretizations based on $\bfH^{\dive}$-conforming elements that independently of how the interface is positioned relative to the computational mesh result in 1) optimal rates of convergence of the approximate velocity and pressure;  2) well-posed linear systems where the condition number of the system matrix scales as it does for fitted finite element discretizations; 3) optimal rates of convergence of the approximate divergence with pointwise divergence-free approximations of solenoidal velocity fields.

In \cite{DanSco12} an unfitted finite element method based on the lowest order Raviart-Thomas space, the $\RT_0\times Q_0$ pair, is developed for the Darcy interface problem. For this element pair the method satisfies property 1) and 3), i.e., has optimal rates of convergence in the $L^2$-norm for velocity, pressure, and divergence, and optimal convergence in the $L^{\infty}$-norm for the divergence.  However, the condition number of the system matrix resulting from the proposed discretization is not controlled and depending on how the interface cuts through the mesh the system matrix may become ill-conditioned. This ill-conditioning is especially pronounced in case of high order elements \cite{puppi2021cut}.

Based on how the interface cuts through the mesh, unfitted finite element methods may produce ill-conditioned linear system matrices reducing the performance of iterative solution techniques~\cite{SauFri13, WadZahKreBer13, puppi2021cut}. To remedy this issue several strategies have been proposed such as: i) adding stabilization terms in the weak form~\cite{Bu10, puppi2021cut, WadZahKreBer13}; ii) using extension operators~\cite{BurHanLar21} to extend the solution from interior elements, that are not cut by the interface, to cut elements;
iii) using agglomeration techniques \cite{JoLa13}; 
iv) using preconditioners, this strategy was for example proposed in \cite{DanSco12} for the $\RT_0\times Q_0$ element pair. See also \cite{GrRe22} and references therein.

 The Cut Finite Element Method (CutFEM) falls into the class of unfitted finite element methods and provides an alternative to standard finite element methods when the construction of a fitted mesh is complicated and expensive.  The basic ideas of CutFEM were first presented in \cite{HaHa02} for an elliptic Partial Differential Equation (PDE). In CutFEM one starts with a regular mesh that covers the computational domain, but does not need to conform to external and/or internal boundaries (interfaces), and appropriate standard finite element spaces.  In this work we consider standard $\bfH^{\dive}$-conforming elements in their pairs such as $\RT_0\times Q_0$, $\BDM_1\times Q_0$, and $\RT_1\times Q_1$~\cite{BofBreFor13}. With these element pairs, standard finite element discretizations satisfy all the desired above mentioned properties in case of a fitted mesh.  
One then defines active unfitted meshes, active finite element spaces, and a weak formulation where interface and boundary conditions are imposed weakly such that the proposed discretization is accurate and robust, independently of the interface position relative to the computational mesh. 
The strategy of adding stabilization terms in the weak form to guarantee a robust discretization independently of the interface position is a popular choice in connection with CutFEM. The reason is that adding ghost penalty stabilization~\cite{Bu10} is relatively easy and has been shown to work well with low as well as high order elements using both continuous and discontinuous elements for bulk problems, see e.g.~\cite{BurHan12, HaLaZa14,  GurStiMas20, PeiFraKreZah22}. 

A cut finite element discretization of the Darcy problem based on Raviart-Thomas spaces and the standard ghost penalty stabilization for both velocity and pressure, in a setting where the mesh is not fitted to the external boundary, has recently been proposed and analyzed in~\cite{puppi2021cut}.  The analysis is for the element pair $\RT_k \times Q_k$, $k\geq 0$. Property 1) and 2) are proven to be satisfied by the numerical scheme~\cite{puppi2021cut}. However, numerical experiments show that the divergence-free property of the Raviart-Thomas element is lost and spurious errors are present in the approximation of the divergence.
 
Recently we proposed a macro-element partitioning of the mesh and could show that, in a cut finite element discretization of a coupled bulk-surface problem, stabilization could be applied more restrictively without losing control of the condition number \cite{LaZa21}.  In this work we show that applying the ghost penalty stabilization only on interior edges of macro-elements reduce and localize errors in the divergence to few macro-elements in the interface or the boundary region.  In order to preserve the optimal approximation of the divergence in the unfitted setting, for the considered element pairs, we propose a modification of the stabilization term for the pressure. Combining the macro-element strategy with the new stabilization terms yields a discretization for which the resulting errors in the divergence are of the same order as when the unstabilized method is used, with the advantage that the linear system is well-posed. The condition number of the system matrix scales with mesh size as it does for fitted finite element discretizations, independently of how the interface cuts through the mesh.  
 
In this paper, we restrict our numerical scheme and the a priori analysis to the case where the mesh is fitted to the external boundary but the interface can cut through the mesh arbitrarily. Essential boundary conditions are imposed strongly while interface conditions are imposed weakly. We derive a priori error estimates for the proposed CutFEM with the element pairs $\RT_k\times Q_k$, $k\geq 0$ and the new stabilization terms. We utilize some results from~\cite{puppi2021cut}.

\textbf{Outline:}
In Section \ref{sec:model} we introduce the interface model problem and its weak formulation. We also set the notation we will use.  In Section \ref{sec:Nummeth} we present two cut finite element discretizations, one based on the standard ghost penalty terms and one based on a new stabilization term for the pressure. We show that the second method, with the new stabilization term, preserves the divergence-free property of the considered elements. In Section~\ref{sec:analysis} we derive a priori error estimates for this method.  In Section \ref{sec:numex} we present numerical experiments and results. Finally, in Section \ref{sec:conclusion} we summarize our findings.
 
\section{The Mathematical Model}\label{sec:model}
We consider a fractured porous medium domain $\Omega \subset \RR^d$, $d=\{2,3\}$, occupied by a fluid wherein the fracture is modelled by the interface $\interf$. We assume $\Omega$ is a bounded convex domain with a polygonal boundary $\partial \Omega$  and $\Gamma$ is a smooth interface that separates $\Omega$ into two disjoint connected subdomains $\Omega_i \subset \Omega$, $i=1,2$, $\Omega=\Omega_1\cup\Omega_2$. When $\interf$ is a closed hypersurface, we take the convention that $\Omega_2$ denotes the region enclosed by $\interf$, with unit normal vector $\bfn=\bfn_1$ of $\interf$ pointing outward from $\Omega_1$ to $\Omega_2$. The notation $\bfn$ is also used for the normal to $\partial \Omega,$ but the context makes it clear which normal we refer to.
We assume the dynamics of the flow is governed by a conservation law and Darcy's law coupled with boundary and interface conditions~\cite{Ber02, MarJafRob05, DanSco12}. Before we state the governing equations we introduce some notation and the problem setup.

\subsection{Notation}
A scalar function $v: \Omega_1 \cup \Omega_2 \rightarrow \RR$ is defined as
\begin{equation}
 v=\begin{cases} v_1,& \text{ in }\Omega_1 \\ v_2,& \text{ in } \Omega_2 \end{cases}, \label{eq:vi}
\end{equation}
with $v_i$ being the restriction of $v$ to $\Omega_i$, $i=1,2$.
The jump operator across the interface $\Gamma$ is defined by
\begin{equation}
\jump{v} = (v_1 - v_2)\restreinta{\Gamma} \label{eq:jump}
\end{equation}
and we denote the average operator by 
\begin{equation}
\{v \}=0.5 (v_1+ v_2)\restreinta{\Gamma}. \label{eq:aver}
\end{equation}

Sometimes we have $v=(v_1, v_2)$ where $v_1$ and $v_2$ overlap in a neighbourhood around $\Gamma$. However, the restriction of $v$ to $\Omega_i$ we define to be $v_i$ and the jump and the average are defined as in equation \eqref{eq:jump} and \eqref{eq:aver}, respectively.

When we write $v \in H^k(\Omega_1 \cup \Omega_2)$ we mean $v=(v_1,v_2)$ where $v_i = v|_{\Omega_i} \in H^k(\Omega_i)$. 
Given $v=(v_1, v_2)$ and $w=(w_1, w_2)$ we will use the notation
\begin{equation} \label{eq:innerp}
(v,w)_{D}=\int_D v(\bfx) w(\bfx)  =\sum_{i=1}^2 \int_{D\cap \Omega_i} v_i(\bfx) w_i(\bfx)  =
\sum_{i=1}^2(v_i,w_i)_{D\cap \Omega_i}
\end{equation}
for the $L^2$-inner product on $D \subset \RR^d$. The norm induced by this inner product will be denoted by $\|\cdot \|_D$. We also have $(v,w)_{\Omega_1 \cup \Omega_2}=\sum_{i=1}^2(v_i,w_i)_{\Omega_i}$. 
We define the pointwise error as 
\begin{equation}
		\| v-w \|_{L^{\infty}({\Omega_1 \cup \Omega_2})} = 
		\max_{i \in \{1,2\}} \max_{\pmb{x}\in{\Omega_i}} |v_i(\bfx)-w_i(\bfx)|.
\end{equation}
For two functional spaces V and W we use the notation $V+W=\{v+w: v\in V, w\in W\}$.

We use similar notations for vector-valued functions and their associated spaces but those functions and spaces are represented using bold letters.

For $\emptyset\neq S\subset \partial D$ smooth enough we define
  \begin{equation}
    \bfH_{0,S}^1(D) := \left\{\pmb{v}:\ \pmb{v} \in  \bfH^1(D), \ \bfv|_S = 0 \right \}.
  \end{equation}
  Recall the functional space 
  \begin{equation}
  \bfH^{\dive}(D) := \left\{\pmb{v}:\ \pmb{v} \in  \mathbf{L}^2(D)  \ \dive \pmb{v} \in L^2(D)\right\}.
    \end{equation}
We also use the spaces 
   \begin{equation}
  \Hdivr(D):= \left\{\pmb{v}:\ \pmb{v} \in  \mathbf{H}^r(D),  \ \dive \pmb{v} \in H^r(D)\right\}, \quad r \in \NN_0
    \end{equation}
   with the norms
  \begin{equation}
  \| \bfv \|_{\Hdivr(D)}^2 = \| \bfv \|_{\bfH^r(D)}^2+\| \dive \bfv \|_{H^r(D)}^2. 
   \end{equation} 
Note that $\bfH^{0,\dive}(D)=\bfH^{\dive}(D)$.

\subsection{The Darcy interface problem} 
Given $\pmb{f}: \Om \rightarrow \RR^d$,  $g: \Om \rightarrow \RR$ , boundary data $\pB: \boundp \rightarrow \RR$, $\uB:  \boundu \rightarrow \RR$ where $\partial \Omega= \boundp  \cup  \boundu$, with  $\boundp  \cap  \boundu = \emptyset$, a prescribed fluid pressure in the fracture $\hat{p}: \interf \rightarrow \RR$, the inverse permeability $\permi$ ($\permi=\bfK^{-1}$ with $\bfK: \Om \rightarrow \RR^d$ being the permeability tensor), and $\permig$ that depends on the permeability tensor in the fracture and the fracture thickness, we seek fluid velocity $\pmb{u}: \Om \rightarrow \RR^d$ and pressure $p: \Om \rightarrow \RR$ satisfying the following Darcy interface problem~\cite{MarJafRob05, DanSco12}: 
\begin{alignat}{2}
  \permi \pmb{u} + \nabla p & = \pmb{f} 
 && \text{ in } \Omega_1\cup\Omega_2, \label{eq:pressure} \\
  \hfill \dive\pmb{u} & =  g
  &&\text{ in } \Omega_1\cup\Omega_2, \label{eq:conservation} \\
  \hfill \jump{p} & =  \permig \{\pmb{u}\cdot\pmb{n}\}
  &&\text{ on } \interf,\label{eq:interfc} \\
  \hfill \{p\} & =  \hat{p}+\xi \permig \jump{\pmb{u}\cdot\pmb{n}}
  &&\text{ on } \interf,\label{eq:interfc2} \\
  \hfill p & =  \pB
  &&\text{ on } \boundp, \label{eq:BCp}\\
  \hfill \pmb{u}\cdot\pmb{n} & =  \uB
  &&\text{ on } \boundu. \label{eq:BCu} 
\end{alignat}
Here, the parameter $\xi\in(0,1/4]$, $\pmb{f} \in \bfL^2({\Omega_1 \cup \Omega_2})$, $g \in L^2({\Omega_1 \cup \Omega_2})$,  $\pB \in H^{1/2}(\boundp)$, $\uB \in H^{-1/2}(\boundu)$, $\hat{p} \in H^{1/2}(\Gamma)$. We assume that the permeability tensor $\bfK$ is diagonal with the smallest diagonal entry being positive and bounded away from zero.

Note that if $\partial \Omega= \boundu$ we have to impose the following compatibility condition on the data: $\int_{\partial \Omega}u_B-\int_\Omega g=0$. To simplify the notation we will hereinafter assume $u_B=0$. Next we formulate a weak formulation of this problem.  

\subsubsection{A weak formulation} \label{sec:weakform}
   Let 
     \begin{align}
     Q & = \left\{q=(q_1,q_2):q_i\in H^1(\Omega_i)\right\}, \\
    \bfV& = \left\{\pmb{v}=(\pmb{v}_1,\pmb{v}_2):\pmb{v}_i\in \bfH^{\dive}(\Omega_i),\ \pmb{v}_i\cdot \pmb{n}|_{\Gamma} \in L^2(\interf), \ \bfv_i\cdot\bfn|_{\boundu}=0 \right\},
     \end{align}
where the traces are defined in the sense of distributions. 
  We now look to find $(\pmb{u},p)\in \bfV\times Q$ such that 
\begin{align}
	a(\pmb{u},\pmb{v}) + b(\pmb{v},p) &= \FF (\pmb{v})  \quad  \forall \pmb{v}\in \bfV, \label{eq:weakdarcy}\\
	b(\pmb{u},q) &= \GG (q) \quad  \forall q\in Q,\label{eq:weakdarcy2}
\end{align}
where
\begin{align}
	a(\pmb{u},\pmb{v}) &:= (\permi \pmb{u},\pmb{v})_{\Omega_1 \cup \Omega_2} +  (\permig \{\pmb{u}\cdot\pmb{n}\},\{\pmb{v}\cdot\pmb{n}\})_{\interf} + (\xi\permig \jump{\pmb{u}\cdot\pmb{n}}, \jump{\pmb{v}\cdot\pmb{n}})_{\interf},
	\label{eq:a}\\
	b(\pmb{u},q) &:= -(\dive \pmb{u},q)_{\Omega_1 \cup \Omega_2}, 
	\label{eq:b} \\
	\FF (\pmb{v}) &:= (\pmb{f},\pmb{v})_{\Omega_1 \cup \Omega_2} 
	- (\pB,\pmb{v}\cdot \pmb{n})_{\boundp}  - (\hat{p}, \jump{\pmb{v}\cdot\pmb{n}})_{\interf}, 
	\label{eq:f} \\
	\GG(q) &:= - (g,q)_{\Omega_1 \cup \Omega_2}. \label{eq:g}
\end{align}

The above weak formulation can be derived using the following integration by parts formula (see e.g. Lemma 2.1.1~\cite{BofBreFor13}):
\begin{equation}  \label{eq:Greens}
 \int_{\Omega_i} \pmb{v} \cdot  \nabla q = - \int_{\Omega_i} \dive \pmb{v} q +  \int_{\partial \Omega_i} \pmb{v} \cdot \pmb{n}_i q, 
\end{equation}
the identity  $\jump{ab}= \jump{a}\{b\}+ \{a\} \jump{b}$, together with the interface and boundary conditions, equations \eqref{eq:interfc}-\eqref{eq:BCu}. Here, $\pmb{v}$ is any function in $\bfV$, $q$ is any function in $Q$, and $\pmb{n}_i$ is the unit normal directed outward from $\Omega_i$. Note that at the interface $\interf$ we have $\pmb{n}=\pmb{n}_1=-\pmb{n}_2$. Using the integration by parts formula on $(\pmb{v},\nabla p)_{\Omega_1 \cup \Omega_2}$ we can rewrite the terms at the interface as
\begin{align} 
   \int_{\interf} p \pmb{v} \cdot \pmb{n}_1+\int_{\interf} p \pmb{v} \cdot \pmb{n}_2 & =   \int_{\interf} \jump{p \pmb{v} \cdot \pmb{n}}  \nonumber \\
   & =  \int_{\interf} \jump{p}\{\pmb{v} \cdot \pmb{n}\}+ \{p\} \jump{\pmb{v} \cdot \pmb{n}} \nonumber \\
   & =  \int_{\interf} \permig \{\pmb{u}\cdot\pmb{n}\}\{\pmb{v} \cdot \pmb{n}\}+ (\hat{p}+\xi \permig \jump{\pmb{u}\cdot\pmb{n}}) \jump{\pmb{v} \cdot \pmb{n}},
\end{align}
and at the boundary we have

\begin{equation} \label{eq:bctreat}
  \int_{\partial \Omega} p\pmb{v} \cdot \pmb{n} = \int_{\boundp} \pB \pmb{v} \cdot \pmb{n}.
\end{equation}

\section{The Finite Element Method}\label{sec:Nummeth}
We now define the meshes, the spaces, and the weak formulations that form our unfitted finite element methods. We present two cut finite element schemes, 
both employing ghost penalty stabilization. 
The first method is based on standard ghost penalty terms both for the pressure and the velocity \cite{BurHanLar21, puppi2021cut} but can use a macro-element partition of the active mesh so stabilization is active only where necessary. The second method contains new mixed stabilization terms related to the $b$-bilinear form, replacing the standard stabilization for the pressure; a macro-element partition can also be utilized here.

For the second method we obtain all the desired properties for all the considered element pairs: convergence, well-conditioned linear systems, and pointwise mass conservation. See Section \ref{sec:numex} for numerical results.

\subsection{Mesh}\label{sec:mesh}
\begin{figure}[h!]
\centering
	\begin{subfigure}[b]{0.45 \textwidth}  	 		
	\centering	
	\includegraphics[scale=0.15]{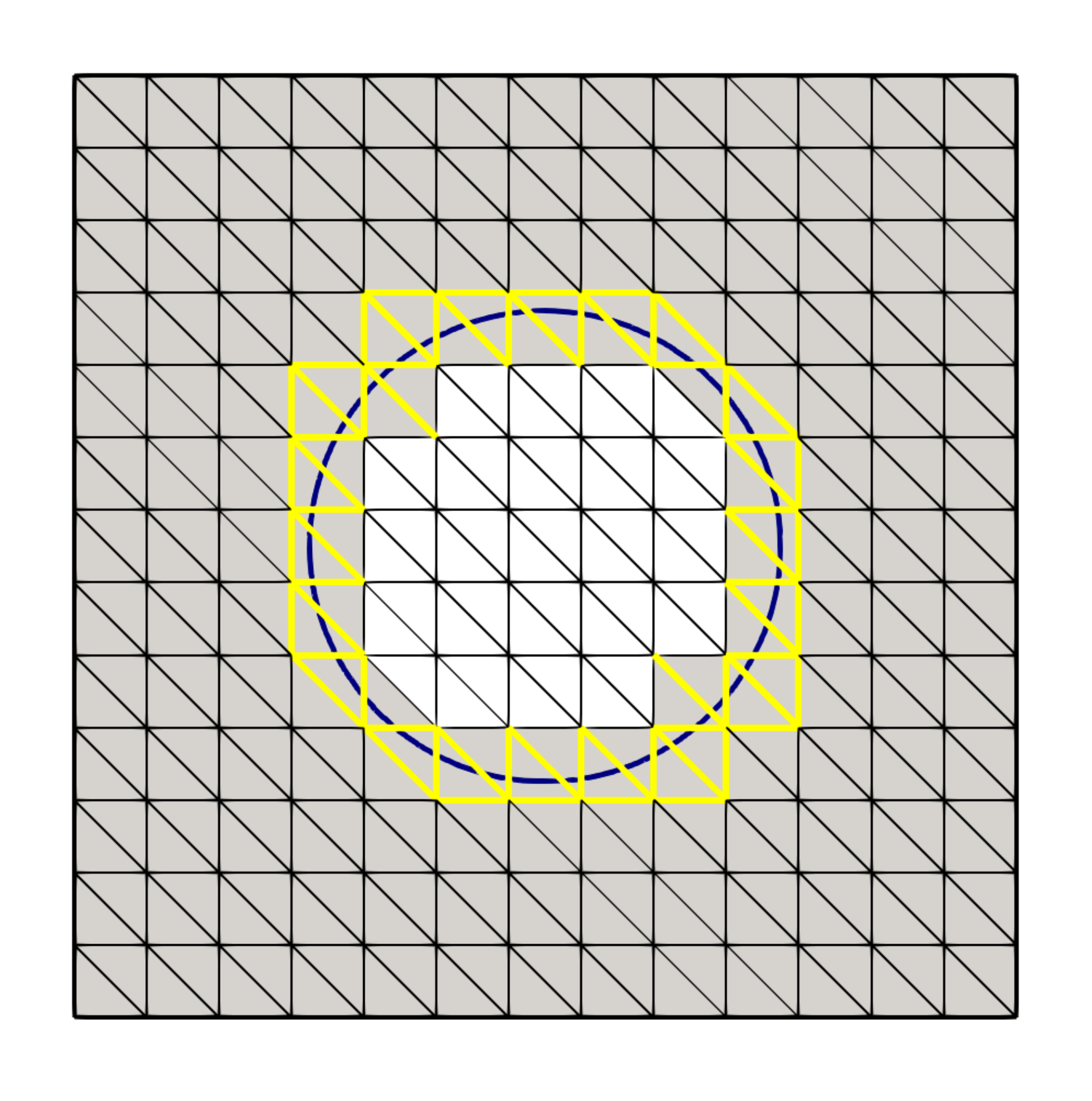}
	\end{subfigure}
	\begin{subfigure}[b]{0.45 \textwidth}  	 	 	
	\centering	
	\includegraphics[scale=0.15]{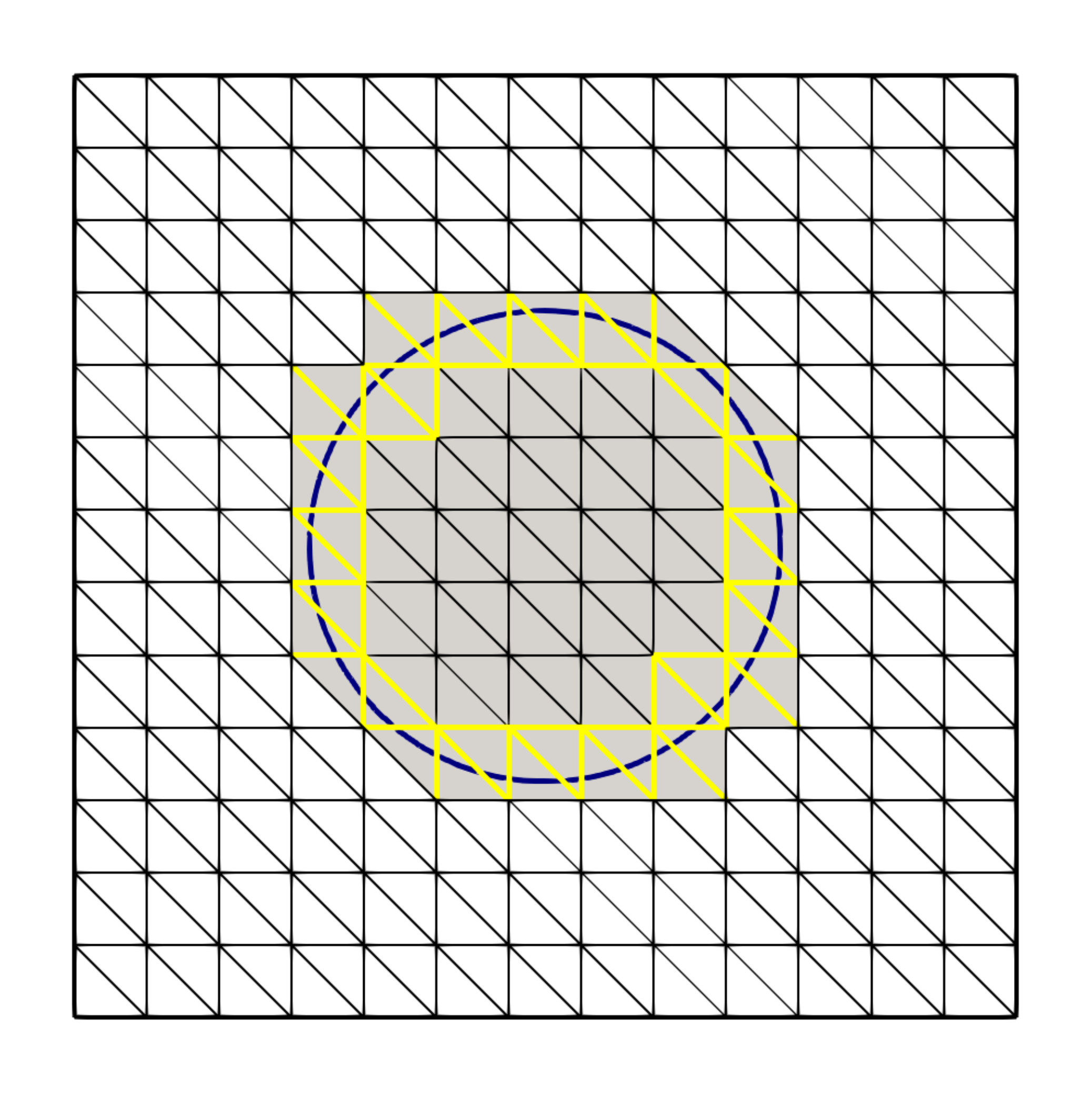}
	\end{subfigure}
      \caption{The grey domain is $\Omega_{\mcT_{h,i}}$ and the yellow faces are faces in $\mcF_{h,i}$. Left:  $i=1$. Right:  $i=2$.  }  
        \label{fig:static active mesh}      	
\end{figure}

Let $\{ \mesh \}_h$ be a quasi-uniform family of simplicial meshes of $\Omega$ with $h\in(0,h_0]$ and $h_0<1$ small. We denote by $h$ the piecewise constant function that on element $T$ is equal to $h_T$, the diameter of element $T\in \mesh$. The mesh $\mesh$, which we refer to as the background mesh, conforms to the fixed polygonal boundary $\partial \Omega$ but does not need to conform to the interface $\interf$. We assume that when the interface cuts an element it intersects the element boundary $\partial T$ exactly twice and each open edge at most once.
Let $\mathcal{G}_h$ be the set of elements intersected by the interface, 
\[
	\mathcal{G}_h = \{T \in \mcT_h : | \Gamma \cap T| > 0 \}
\]
and $\mcF_h$ be the set of faces (edges) in $\mathcal{G}_h$. 
We define active meshes $\mcT_{h,i}$ and domains $\Omega_{\mcT_{h,i}}$ associated to the subdomains $\Omega_i$ by
\begin{align}\label{eq:activemd}
  \mcT_{h,i} = \{T \in \mcT_h : T \cap \Omega_i \neq \emptyset \}, \quad 
  \Omega_{\mcT_{h,i}}:=\bigcup_{T\in\mcT_{h,i}}T,  \quad i=1,2. 
\end{align}
We denote by $\mcF_{h,i}$ the set of faces in $\mcF_h$ that are shared by two elements of $\mcT_{h,i}$. These sets are illustrated in Figure \ref{fig:static active mesh}.\\

\subsection{Spaces}
Define the following finite element spaces on the active meshes $\mcT_{h,i}$:

\begin{align} \label{eq:Vi}
	\Vki & =\left\{\pmb{v}_h\in \bfH^{\dive}(\Omega_{\mcT_{h,i}}): \ \pmb{v}_h\lvert_T\ \in \bfV_k(T), \ \forall T\in \mcT_{h,i},\ \bfv_h\cdot\bfn|_{\boundu} = 0 \right\},  
\end{align}

\begin{align}
  Q_{\kp}(\mcT_{h,i}) & =\bigoplus_{T \in \mcT_{h,i}} Q_{\kp}(T),
\end{align}
where $Q_{\kp}(T)$ is the space of polynomials in $T$ of degree less or equal to $\kp\geq 0$ and $i=1,2$. Here $\bfV_k=\RT_{k}$ or $\bfV_k=\BDM_k$. Having chosen $\bfV_k$ then $\kp$ is implicitly chosen: when $\bfV_k=\RT_{k}$, the Raviart-Thomas space, we take $\kp=k\geq 0$ and when $\bfV_k=\BDM_k$, the Brezzi-Douglas-Marini space, we take $\kp=k-1 \geq 0$. For the definition of these spaces we refer to for example Chapter III.3 in \cite{BrFo12}.

Finally, let
\begin{align*}
 \Vk&=\left\{\pmb{v}_h = (\pmb{v}_{h,1},\pmb{v}_{h,2}) : \ \pmb{v}_{h,i} \in \Vki, \ i=1,2 \right\} , \\
  \Qk &= \left\{q_h = (q_{h,1},q_{h,2}): \ q_{h,i} \in \Qki, \ i=1,2 \right\}.
\end{align*}
If $\boundp=\emptyset$ so that $\partial \Omega = \boundu$ we define 
\begin{align}\label{eq:Qk_pure}
  \Qk &= \left\{q_h = (q_{h,1},q_{h,2}): \ q_{h,i} \in \Qki \ \text{ and } \int_{\Omega_1 \cup \Omega_2} q_h=0  \right\}.
\end{align}


\subsection{Weak forms}
We introduce an unstabilized unfitted discretization as a starting point for the stabilized schemes:  
Find $(\pmb{u}_h,p_h)\in \Vk \times \Qk$ such that
\begin{align}\label{eq:unstabdarcy}
	a(\pmb{u}_h,\pmb{v}_h) + b(\pmb{v}_h,p_h) &= \FF (\pmb{v}_h) \quad \forall \pmb{v}_h\in \Vk, \\
	b(\pmb{u}_h,q_h) &= \GG (q_h) \quad \forall q_h\in\Qk.
\end{align}
Note that in this work the mesh is fitted to $\partial\Omega$ and we impose the essential boundary conditions in the functional space $\Vk$. To simplify the notation we have assumed homogeneous boundary conditions. In case of non-homogeneous boundary conditions one needs to change the trial space and impose the boundary condition $u_B$ but the test space remains the same.

We now present two cut finite element discretizations. Stabilization terms are added in the variational formulation to provide stability and control of the condition number independent of how the background mesh is cut by the interface. In the first method we add standard ghost penalty terms, both for velocity and pressure, while in the second method we introduce alternative ghost penalty terms. The second method will have all the desired properties.

\subsubsection{Method 1: Standard ghost penalty terms} \label{sec:methodstab}
The finite element method reads: Find $(\pmb{u}_h,p_h)\in \Vk\times \Qk$ such that

\begin{equation}\label{eq:discretedarcy}
    a(\pmb{u}_h,\pmb{v}_h) + b(\pmb{v}_h,p_h) - b(\pmb{u}_h,q_h) + s_{\pmb{u}}(\pmb{u}_h,\pmb{v}_h) +  s_{p}(p_h,q_h) = \FF (\pmb{v}_h)-\GG(q_h),   
\end{equation}

for all $(\pmb{v}_h,q_h)\in \Vk\times \Qk$. 
Here $s_{\pmb{u}}$ and $s_p$ are ghost penalty stabilization terms \cite{Bu10}, connecting elements with a small intersection with $\Omega_i$ to elements with a large intersection,
\begin{align}
s_{\pmb{u}}(\pmb{u}_h,\pmb{v}_h) &= \sum_{i=1}^2 \sum_{F \in \mathcal{F}_{h,i}}\sum_{j=0}^{\hat{k}+1}
\tau_{\bfu} h^{2j+1}  (\jump{D^j_{\bfn_F} \pmb{u}_{h,i} }, \jump{D^j_{\bfn_F} \pmb{v}_{h,i}})_{F}, \label{eq:full stab-u method1} \\
s_p(p_h,q_h) &= \sum_{i=1}^2 \sum_{F \in \mathcal{F}_{h,i}}\sum_{j=0}^{\hat{k}}
\tau_p h^{2j+\gamma}  (\jump{D^j p_{h,i} }, \jump{D^j q_{h,i}})_{F}. \label{eq:full stab-p method1}
\end{align}
Here $\jump{D^j p_{h,i}}$ denotes the jump in the derivative of order $j$ across the face $F$, $\jump{D^0_{\bfn_F} \bfv_{h,i}} = \jump{\bfv_{h,i}}$, and 
\begin{equation*}
  (\jump{D^j_{\bfn_F} \pmb{u}_{h,i} }, \jump{D^j_{\bfn_F} \pmb{v}_{h,i}})_{F}= \sum_{l=1}^d(\jump{D^j_{\bfn_F} \pmb{u}_{h,i}^l }, \jump{D^j_{\bfn_F} \pmb{v}_{h,i}^l})_{F},
\end{equation*}
where $D^j_{\bfn_F}$ denotes the normal derivative of order $j$. 
The stabilization parameters $\tau_{\pmb{u}}>0$, and $\tau_{p}>0$ may vary with $j$ and $i$, and $\gamma \in [-1 ,1]$ is a parameter. We choose $\gamma=1$ so that the scheme has both optimal approximation properties and the condition number of the associated linear system matrix scales as for standard fitted finite element discretizations, see e.g. \cite{arnold_precond}.

The stabilization terms extend the approximated velocity and pressure from $\Omega_i$ to the entire active domain $\Omega_{\mcT_{h,i}}$, $i=1,2$. We state a result on the control of the $L^2$-norm of functions in  $\Qki$.
Here and after we use the notation $x_1\lesssim x_2$ if and only if $x_1\leq C x_2$ for some constant $C$ that is independent of the mesh parameter $h$ and how the interface cuts through the mesh.
\begin{lemma}\label{lem:sp_ineq}
  For any $q_{h,i} \in \Qki$, 
  \[ \|q_{h,i} \|^2_{\Omega_{\mcT_{h,i}}} \lesssim \|q_{h,i}\|^2_{\Omega_i} + s_p(q_{h,i},q_{h,i}) . \]
\end{lemma}
The proof follows from the arguments in the proof of Lemma 3.8 in \cite{HaLaZa14} or Lemma 5.1 in  \cite{MaLaLoRo14}. A corresponding result holds for functions in $\Vki$ using the stabilization term $s_{\bfu}$. 

\subsubsection{Method 2: New ghost penalty terms}\label{sec:method2}
The finite element method now reads: Find $(\pmb{u}_h,p_h)\in \Vk \times \Qk$ such that
\begin{align}
	\Af(\pmb{u}_h,\pmb{v}_h) + \Bf(\pmb{v}_h,p_h) &= \FF (\pmb{v}_h) \quad \forall \pmb{v}_h\in \Vk,\label{eq:discretedarcy2} \\
	\Bf(\pmb{u}_h,q_h) &= \GG (q_h) \quad \forall q_h\in \Qk.\label{eq:discretedarcy2b}
\end{align}
Here  
\begin{align}
	&\Af(\pmb{u}_h,\bfv_h) := a(\pmb{u}_h,\bfv_h) + s_{\pmb{u}}(\pmb{u}_h,\bfv_h) \nonumber \\
	& \ =(\permi \pmb{u}_h,\bfv_h)_{\Omega_1 \cup \Omega_2} +  (\permig \{\pmb{u}_h\cdot\pmb{n}\},\{\bfv_h\cdot\pmb{n}\})_{\interf} + (\xi\permig \jump{\pmb{u}_h\cdot\pmb{n}}, \jump{\bfv_h\cdot\pmb{n}})_{\interf}+ s_{\bfu}(\pmb{u}_h,\bfv_h), \\
	& \Bf(\pmb{u}_h,q_h) :=b(\pmb{u}_h,q_h)-s_b(\pmb{u}_h,q_h) = -(\dive \pmb{u}_h,q_h)_{\Omega_1 \cup \Omega_2} - s_b(\pmb{u}_h,q_h), \label{eq: Bf0} 
\end{align}
with $s_{\pmb{u}}(\pmb{u}_h,\bfv_h)$ as in Method 1, see equation \eqref{eq:full stab-u method1}, and 
\begin{align}
s_b(\pmb{u}_h,q_h) &= \sum_{i=1}^2 \sum_{F \in \mathcal{F}_{h,i}}\sum_{j=0}^{\hat{k}}
\tau_b h^{2j+1}  ( [D^j (\dive \pmb{u}_{h,i})], [D^j q_{h,i}] )_{F}. \label{eq:full stab-p method2}
\end{align}
The stabilization parameter $\tau_{b}>0$ can vary with $j$ and $i$.

The linear system of equations associated with \eqref{eq:discretedarcy2}-\eqref{eq:discretedarcy2b} can be written as
\begin{equation*}
	\begin{bmatrix}
		\mathbf{A} & \mathbf{B}^T \\
		\mathbf{B} & \mathbf{0}
	\end{bmatrix} 
	\begin{bmatrix}
		\mathbf{\hat{u}} \\
		\mathbf{\hat{p}} 
	\end{bmatrix} 
	=\begin{bmatrix}
		\mathbf{\hat{f}} \\
		\mathbf{\hat{g}} 
	\end{bmatrix}, 
\end{equation*}
with unknown degree of freedom (DOF) values $\begin{bmatrix}
	\mathbf{\hat{u}} \\
	\mathbf{\hat{p}} 
\end{bmatrix}$ and data $\begin{bmatrix}
\mathbf{\hat{f}} \\ 
\mathbf{\hat{g}}  
\end{bmatrix}$. 
Note that $s_b(\bfv_h,p_h)=s_p(\dive \bfv_h, p_h)$ and with this mixed ghost penalty term, Method 2 fits into the framework of symmetric saddle-point problems. This is not the case for Method 1 which utilizes the standard ghost penalty term $s_p(p_h,q_h)$. The stabilization term $s_b(\bfv_h,p_h)$ yields control of the approximated pressure in the active mesh as $s_p$ does but is chosen so that the divergence-free property of $\bfH^{\dive}$-conforming elements is not lost. Next we show that Method 2 produces pointwise divergence-free approximations of solenoidal velocity fields, as fitted standard FEM does with the considered element pairs.

\begin{thm} \label{thm:divfree}(\textbf{The divergence-free property of Method 2}) Assume $g=0$ and $\bfu_h=(\pmb{u}_{h,1},\pmb{u}_{h,2}) \in\Vk$ satisfies \eqref{eq:discretedarcy2}-\eqref{eq:discretedarcy2b}. Then $\dive\bfu_{h} =0$.
\end{thm} 
\begin{proof}
We have that $\bfu_h$ satisfies 
	\[0 = \Bf(\pmb{u}_h,q_h) = -\int_{\Omega_1 \cup \Omega_2} \dive\bfu_h q_h - \sum_{i=1}^2 \sum_{F \in \mathcal{F}_{h,i}}\sum_{j=0}^{\hat{k}}
	\tau_b h^{2j+1}  ( \jump{D^j (\dive \pmb{u}_{h,i})}, \jump{D^j q_{h,i}} )_{F},  \quad \forall  q_h\in\Qk.  \]
	We may choose $q_h=(q_{h,1},q_{h,2})$ such that $q_{h,i}=-\dive\bfu_{h,i}$ since $\dive \left( \Vki \right) \subseteq \Qki$. Then we have 
	\begin{equation}
		0=\|\dive\pmb{u}_{h,i} \|_{{\Omega_i}}^2+ s_p(\dive\pmb{u}_{h,i}, \dive\pmb{u}_{h,i}) 
		\gtrsim \|\dive\pmb{u}_{h,i} \|_{{\Omdh}}^2\geq 0,
	\end{equation}
	where we used that for any $q_{h,i} \in\Qki$, $\|q_{h,i} \|_{\Omega_i}^2 + s_p(q_{h,i}, q_{h,i}) \gtrsim \| q_{h,i} \|^2_{\Omdh}$, see Lemma \ref{lem:sp_ineq}. We can thus immediately conclude that $\dive\bfu_{h,i} =0$ in $\Omdh$. 
\end{proof}


\subsection{Macro-element stabilization} \label{sec:macro-elements}
Stenberg introduced a macro-element technique, see e.g. \cite{Ste84}, based on decomposing the mesh into disjoint macro-elements. We use this idea, but in a different setting, to create elements with a large $\Omega_i$-intersection. As introduced in \cite{LaZa21}, the idea is that if each macro-element has a large $\Omega_i$-intersection stabilization is only needed on internal faces of macro-elements and never on the interface between two macro-elements. We thus stabilize on fewer faces and hence get fewer non-zero entries in the system matrix. In \cite{LaZa21} it was also illustrated that with the macro-element stabilization the error from the discretization is not sensitive to large stabilization parameters ($\tau_{\bfu}$, $\tau_p$, $\tau_b$).

\begin{figure}[h!]
\centering
	\begin{subfigure}[b]{0.4 \textwidth}  	 	 	
	\centering	
	\includegraphics[scale=0.165]{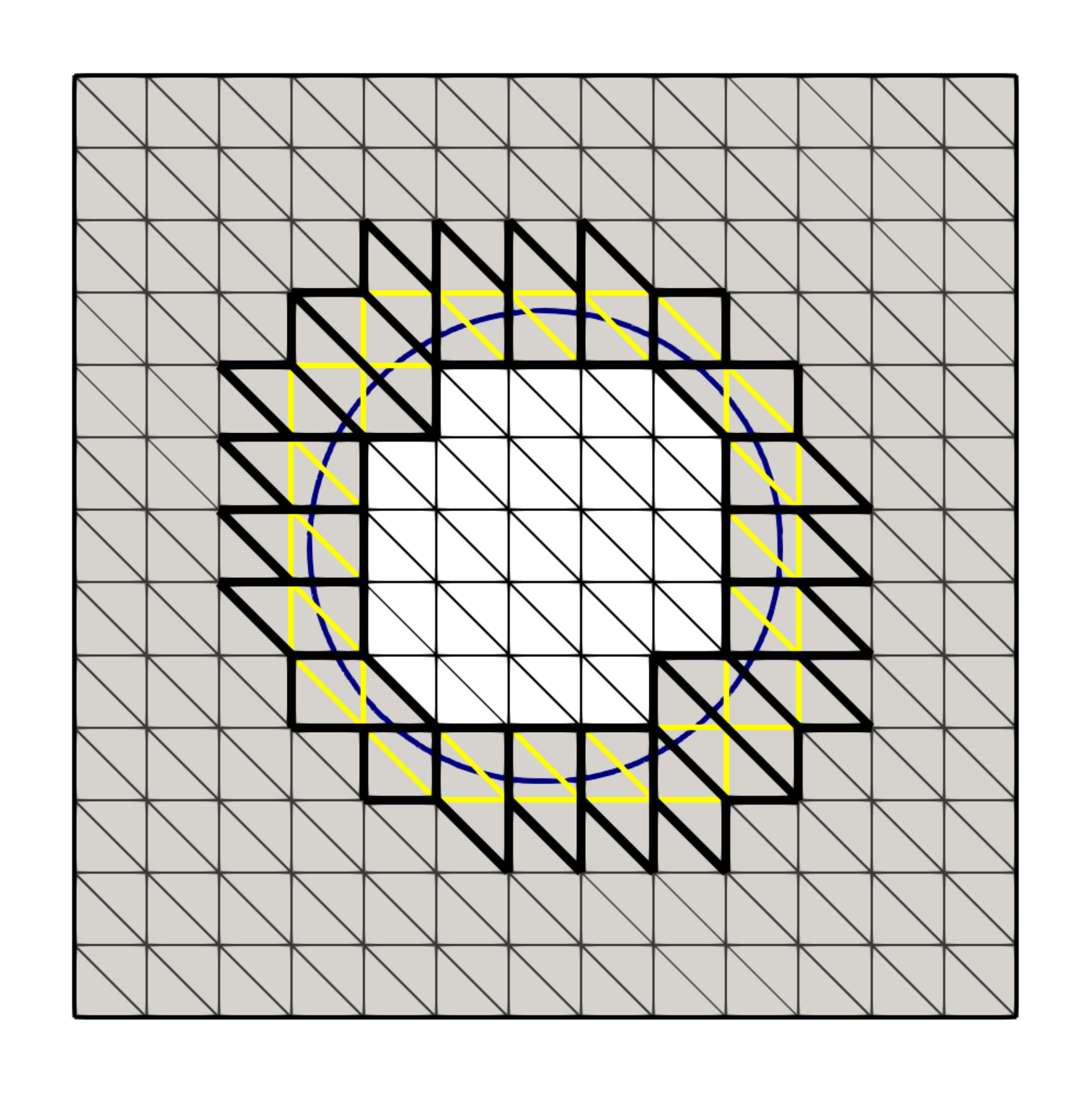}
	\end{subfigure}
	\begin{subfigure}[b]{0.4 \textwidth}  	 	 	
	\centering	
	\includegraphics[scale=0.165]{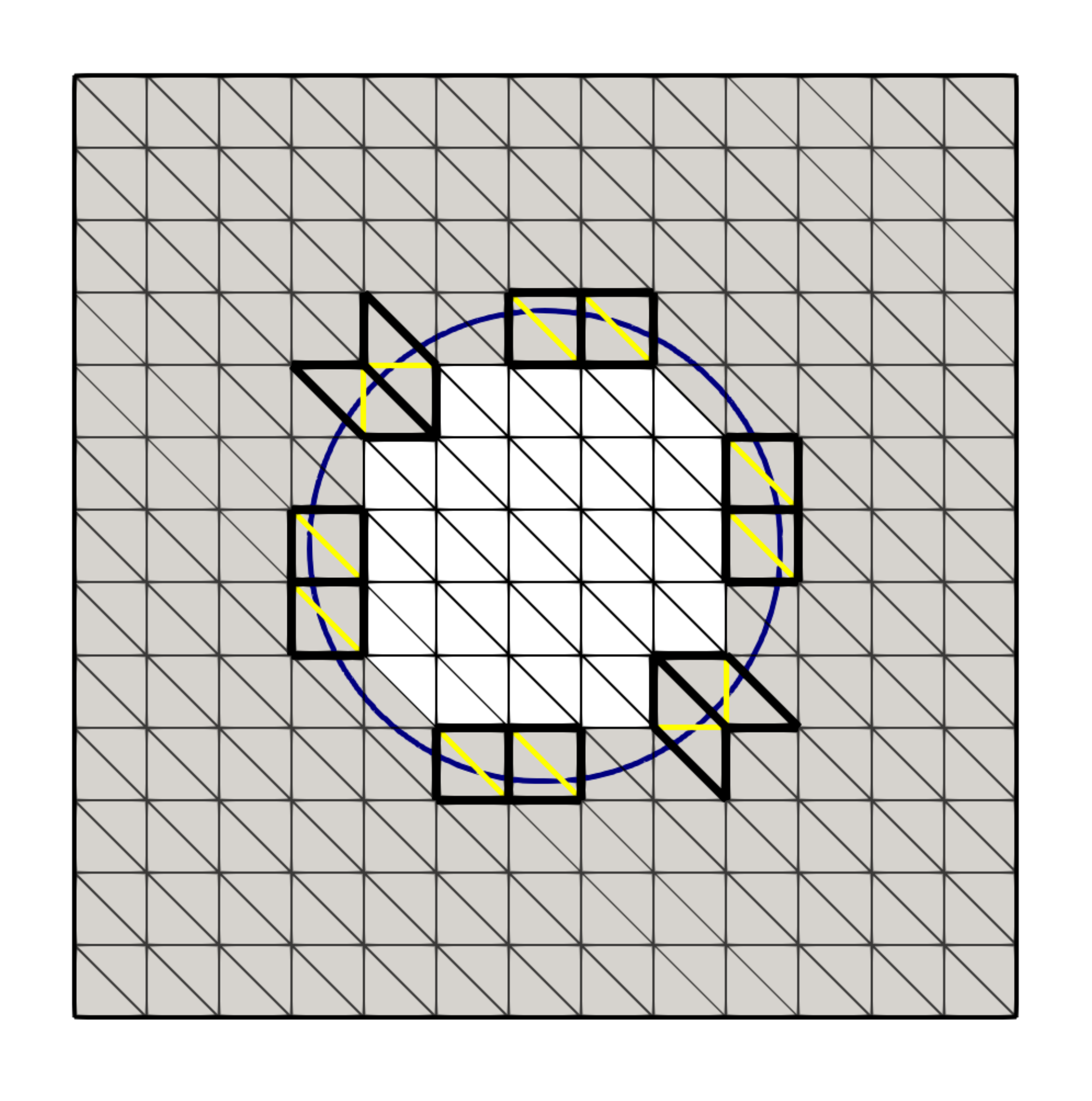}
	\end{subfigure}
        \caption{Decomposition of the triangulation $\mcT_{h,1}$ into macro-elements for different $\delta_1$ (see \eqref{eq:largeel}).
   The thick black lines illustrate the boundary of macro-elements consisting of more than one element. The yellow lines are interior faces of the macro-elements. Left: $\delta_1 = 1$. Right: $\delta_1 = 0.2$.}
         \label{fig:macro_mesh} 
\end{figure}

We start by classifying each element in the active mesh $\mcT_{h,i}$ as either having a large intersection with the domain $\Omega_i$ or a small intersection. An element $T \in \mcT_{h,i}$ has a large $\Omega_i$-intersection if
\begin{equation}\label{eq:largeel}
  \delta_i  \leq \frac{|T \cap \Omega_i|}{|T|},
\end{equation}
where $0<\delta_i\leq 1$ is a constant that is independent of the element and the mesh parameter. We collect all elements with a large intersection in 
\begin{equation}
  \mcTL=\left\{T \in \mcT_{h,i} : |T \cap \Omega_i| \geq \delta_i |T| \right\}.
\end{equation}
Using this classification we create a macro-element partition $\mcM_{h,i}$ of $\Omega_{\mcT_{h,i}}$ following \cite{LaZa21}:
\begin{itemize}
\item To each $T_L\in\mcTL$ we associate a macro-element mesh $\mcT_{h,i}(T_L)$ containing $T_L$ and possibly adjacent elements that are in $\mcT_{h,i}\setminus \mcTL$, i.e., elements classified as having a small intersection with $\Omega_i$ and connected to $T_L$ via a bounded number of internal faces. 
\item Each element $T \in \mcT_{h,i}$ belongs to precisely one macro-element mesh $\mcT_{h,i}(T_L)$.
\item Each macro-element $M_L \in \mcM_{h,i}$ is the union of elements in $\mcT_{h,i}(T_L)$, i.e., 
\begin{equation}
M_L = \bigcup_{T \in \mcT_{h,i}(T_L)} T.
\end{equation}
\end{itemize}
We denote by $\mcFM$ the set consisting of interior faces of $M_L\in\mcM_{h,i}$. Note that $\mcFM$ is empty when $T_L$ is the only element in $\mcTT$.
See Figure~\ref{fig:macro_mesh} for an illustration of a macro-element partitioning. Faces in $\mcFM$ are shown in yellow. 

Stabilization is then applied only on internal faces of macro-elements and never on faces shared by neighbouring macro-elements. The stabilization corresponding to \eqref{eq:full stab-u method1}, \eqref{eq:full stab-p method1}, and \eqref{eq:full stab-p method2}  becomes
\begin{align}
  s_{\pmb{u}}(\pmb{u}_h,\pmb{v}_h) &= \sum_{i=1}^2 \sum_{M_L \in \mcM_{h,i}} \sum_{F \in \mcFM} \sum_{j=0}^{\kp+1} \tau_{\bfu} h^{2j + 1} (\jump{D^j_{\bfn_F} \pmb{u}_{h,i} }, \jump{D^j_{\bfn_F} \pmb{v}_{h,i}})_{F}, \\
  s_{p}(p_h,q_h) &= \sum_{i=1}^2 \sum_{M_L \in \mcM_{h,i}}  \sum_{F \in \mcFM} \sum_{j=0}^{\kp} \tau_{p} h^{2j + 1} (\jump{D^j p_{h,i}}, \jump{D^j q_{h,i}})_{F}, \\
  s_{b}(\pmb{u}_h,q_h) &= \sum_{i=1}^2 \sum_{M_L \in \mcM_{h,i}} \sum_{F \in \mcFM}  \sum_{j=0}^{\kp} \tau_{b} h^{2j + 1} ( \jump{D^j (\dive \pmb{u}_{h,i})}, \jump{D^j q_{h,i}} )_{F}. \label{eq:stab-b}
\end{align}
Note that with the macro-element stabilization we apply stabilization more restrictively, compare the number of yellow faces in Figure~\ref{fig:macro_mesh} illustrating faces involved in the macro-element stabilization with Figure~\ref{fig:static active mesh} which illustrates faces involved in full stabilization.  From \cite[Lemma 3.2]{LaZa21} it follows that Lemma~\ref{lem:sp_ineq} also holds with the macro-element stabilization. We follow Algorithm 1 in \cite{LaZa21} when we construct the macro-element partition. 

\section{Analysis of the Method}\label{sec:analysis}
Proofs for well-posedness of the unstabilized unfitted method \eqref{eq:unstabdarcy} exist only for the lowest order pair $\RT_0\times Q_0$ \cite{DanSco12}, to the best of our knowledge. Despite this, using a direct solver, the unstabilized method succeeds in approximating the exact solution in all our numerical experiments with optimal convergence order. However, the resulting linear systems could be severely ill-conditioned in particular in case of high order elements, we illustrate this in the next Section. The analysis of Method 1 would to a large extent follow from the work in~\cite{puppi2021cut}. However, this method does not preserve the divergence-free property of the $\bfH^{\dive}$-conforming elements. We propose to use Method 2. Thus, we only analyze Method 2 and throughout this section we assume $\bfV_k=\RT_k$ and hence $\hat{k}=k$.

We start by defining the following mesh dependent norms: 
\begin{align}
\| q \|^2_h &:= 
	\sum_{i\in\{1,2\}}\sum_{T\in\mcT_{h,i}} \| q_{i} \|^2_T = \| q \|^2_{\Omdh}, \\
    \tn \bfv \tn^2_\Af &:= \Af(\bfv,\bfv) \nonumber \\
    &= \|\permi^{1/2}\bfv\|^2_{\Omega_1\cup\Omega_2} +\|\permig^{1/2}\{\bfv\cdot\pmb{n}\}\|^2_{\interf} + \xi \|\permig^{1/2}\jump{\bfv\cdot\pmb{n}}\|^2_{\interf} + s_{\bfu}(\bfv,\bfv),  \label{eq:anorm}\\
	\tn \bfv \tn^2_h &:= \tn \bfv \tn_\Af^2 + \| \dive \bfv\|_{\Omdh}^2=\tn \bfv \tn_\Af^2 + \| \dive \bfv\|_{h}^2.\label{eq:vhnorm}
\end{align}

Next we state some properties of standard interpolants for our spaces that we will use.  

Let $\projpi:L^2(\Omdhi)\to \Qki$ be the orthogonal projection on $\Qki$ and $\interpvi : \bfH^{\dive}(\Omdhi)\cap \prod_{T\in\mcT_{h,i}} \bfH^r(T) \to \Vki$, ($r\geq 1$), be the global interpolation operator associated to $\Vki$, see e.g. \cite[(2.5.10)-(2.5.27)]{BofBreFor13} and \cite[Chapter 3]{Gatica14}. The following commuting diagram holds (see \cite[diagram (2.5.27)]{BofBreFor13} or \cite[Lemma 3.7]{Gatica14}):
\begin{equation}\label{eq:diagram}
	\begin{tikzcd}
		\bfH^{\dive}(\Omega_{\mcT_{h,i}})\cap \prod_{T\in\mcT_{h,i}} \bfH^{r}(T) \arrow{r}{\dive} \arrow[swap]{d}{\pi_{h,i}^k} & L^2(\Omega_{\mcT_{h,i}}) \arrow{d}{\projpi} \\%
		\Vki \arrow{r}{\dive} & \Qki
	\end{tikzcd}
\end{equation}


The following local interpolation error estimates hold:
\begin{lemma}\textbf{(Local interpolation error estimates)}\label{lem:loc_interp}
	Fix $T\in\mcT_{h}$. Let $\bfu\in\bfH^{k+1}(T)$, $\dive\bfu\in H^{k+1}(T)$, then for $0\leq m\leq \ell+1$ and $0\leq \ell \leq k$,
	\begin{align*}
	| \bfu-\interpvi \bfu |_{m,T} &
        \lesssim h^{\ell+1-m}|\bfu|_{\ell+1,T}, \\
	| \dive\bfu-\dive\interpvi \bfu |_{m,T} &
        \lesssim h^{\ell+1-m}|\dive\bfu|_{\ell+1,T}.
	\end{align*} 
\end{lemma} 
For a proof see Lemma 3.16 and 3.17 of \cite{Gatica14}.

We will also make use of the following elementwise trace inequalities. For $T\in \mesh$  we have the trace inequalities
\begin{equation} \label{eq:tracestand}
\| v \|^2_{L^2(\partial T)} \lesssim h^{-1} \| v \|_{L^2(T)}^2 + h \| \nabla v \|^2_{L^2(T)} \qquad v\in H^1(T), 
\end{equation} 
\begin{equation} \label{eq:trace}
\| v \|^2_{L^2(T\cap \Gamma)} \lesssim h^{-1} \| v \|_{L^2(T)}^2 + h \| \nabla v \|^2_{L^2(T)} \qquad v\in H^1(T),
\end{equation} 
where the first inequality is a standard trace inequality, see, e.g.~\cite{BreSco}, and the second inequality is proven in \cite{HaHaLa04} and the constant 
is independent of how $\Gamma$ intersects $T$ and of $h$.

\subsection{Stability} \label{sec:stability}
Before we state and prove an inf-sup result for $\Bf$ we state two Lemmas we will need.    
\begin{lemma}\textbf{(The Poincaré inequality)}\label{lem:Poincare}
	Let $D$ be a bounded Lipschitz domain in $\RR^d$
        with $\emptyset\neq S\subset \partial D$
        and let $\bfv\in\bfH_{0,S}^1(D)$. Then 
	\begin{align*}
		\|\bfv\|_{D} \leq C_P \|\nabla\bfv\|_{D},
	\end{align*} 
	where $C_P=2$diam$(D)$.
\end{lemma}
See \cite[Theorem 7.91]{Salsa2016}.

 
\begin{lemma}\label{lem:divergence_Poincare}
 Let $D$ be a bounded convex polygonal domain in $\RR^d$ or a bounded open subset of $\RR^d$ with a smooth boundary $\partial D$.
        Then there exist a positive constant C such that for all $q\in L^2(D)$ there is a $\bfv \in \bfH^1(D)$ satisfying 
	\[\dive\bfv = q,\ \|\bfv\|_{\bfH^1(D)}\leq C \|q\|_{D}. \]
        If $q$ satisfies $\int_D q=0$ then $\bfv \in \bfH_{0}^1(D)$.
\end{lemma} 
See \cite[Lemma 11.2.3]{BreSco}. We also refer to \cite{GirRav86}.


Let $U_{\varepsilonz}(\Gamma)$ be an open tubular neighborhood of $\Gamma$ with thickness $\varepsilonz>0$ such that for each $\bfx \in U_{\varepsilonz}(\Gamma)$ there is a uniquely determined $p(\bfx) \in \Gamma$. We state and prove an inf-sup result under the assumption that $\Gamma$ is smooth and simply connected and such that $U_{\varepsilon}(\Gamma)$ with $\varepsilon=ch \leq \varepsilon_0$ satisfies the conditions for the domain $D$ of Lemma \ref{lem:divergence_Poincare}. We also assume that $h$ is sufficiently small so that $\bigcup_{T \in \mathcal{G}_h } T \subset U_{\varepsilon}(\Gamma)$  and hence $\Omdhi \subset \Omeps=\Omega_i \cup U_{\varepsilon}(\Gamma)$. Note that $\varepsilonz$ is a constant that depends on the curvature of the surface $\Gamma$. For a surface with high curvature the constant $\varepsilonz$ is small.

\begin{prop}(\textbf{The inf-sup condition}) \label{prop:infsup}
For every $q_h\in \Qk$ there exists $\bfv_h \in \Vk$ such that 
\begin{align}
\Bf(\bfv_h,q_h) \gtrsim \| q_{h} \|^2_h, \quad \tn \bfv_h \tn_h \lesssim \| q_h\|_h.
\label{eq:inf sup proposition} 
\end{align}
\end{prop}

\begin{proof}
Fix $q_{h}=(q_{h,1},q_{h,2}) \in \Qk$. Let
\begin{equation}
  \hat{q}=
  \left \{
  \begin{tabular}{ll}
  $q_{h,1}$ & \textrm{in $\Omega_1$,} \\
  $q_{h,2}$ & \textrm{in $\Omega_2$.}
  \end{tabular}
  \right.
\end{equation}
Then, for $\hat{q} \in L^2(\Omega)$, let $\bfv_p \in \bfH^{1}(\Omega)$
be the function from Lemma \ref{lem:divergence_Poincare} satisfying 
\begin{align} 
&\dive \bfv_p =-\hat{q} \quad \textrm{ in $\Omega$,}   \label{eq:vp} 
\\
&\|\bfv_p\|_{\bfH^1(\Omega)} \lesssim \|\hat{q}\|_{\Omega}. 
\label{eq: vpnorm}
\end{align} 
Note that in case $\boundp=\emptyset$ we have $\int_{\Omega_1 \cup \Omega_2} q_h=0$, (i.e., $\int_{\Omega} \hat{q}=0$) then according to Lemma \ref{lem:divergence_Poincare} $\bfv_p \in \bfH^{1}_{0}(\Omega)$. This is consistent with integration by parts
  \[0=\int_{\partial\Omega}\bfv_p\cdot\bfn=\int_{\Omega}\dive\bfv_p = -\int_{\Omega} \hat{q}=0.\]  
  Further for each $i=1,2$ extend $q_{h,i}$ defined in $\Omdhi$ by zero and denote this extension by $q_{h,i}^e$.  
  Following the proof of Lemma 4.7 in \cite{DanSco12} we can construct a  $\bfv_i^e \in \bfH^{1}_{0,S_i}(\intfneig)$ with $\emptyset\neq S_i \subset\partial \intfneig$ ($S_i$ belongs to the part of the boundary that is outside $\Omdhi$)  such that 
\begin{align}
&\dive \bfv_i^e =-q_{h,i}^e \quad \textrm{ in $\intfneig$,}   \label{eq:vi} 
\\ 
&\|\bfv_i^e\|_{\bfH^1(\intfneig)} \lesssim \|q_{h,i}^e\|_{\intfneig}. \label{eq:vinorm}
\end{align} 
The idea is that we can construct $q_{0,i}=q_{h,i}^e- \bar{q}_i \dive \bfh_i$ with $\bar{q}_i=\int_{\intfneig}q_{h,i}^e$ and $\bfh_i \in \bfH^{1}_{0,S_i}(\intfneig)$. Note that $\int_{\partial \intfneig} \bfh_i \cdot \bfn\neq 0$ and we can require $\int_{\intfneig} \dive \bfh_i =1$, thus $\int_{\intfneig} q_{0,i}=0$. Then from Lemma \ref{lem:divergence_Poincare} we know that for $q_{0,i}$ there is a $\bfv_{0,i} \in \bfH^{1}_{0}(\intfneig)$ and $\bfv_i^e=\bfv_{0,i}+\bar{q}_i\bfh_i$ is our desired function.

Let $\bfv_{h}=(\bfv_{h,1}, \bfv_{h,2}) \in \Vk$ be the Raviart-Thomas interpolant of $\bfv_p$ on the uncut elements in $\Omega_{h,i}$ and of $\bfv_i^e$ on the cut elements:   
\begin{equation}\label{eq:defvhi}
  \bfv_{h,i} = \interpvi \bfv_i
  \quad \textrm{with }
 \bfv_i =  \left \{
  \begin{tabular}{ll}
   $\bfv_p$  & \textrm{for $T \in \mcT_{h,i} \setminus \mathcal{G}_h$,} \\
   $\bfv_i^e$ & \textrm{for $T \in \mathcal{G}_h$.}
  \end{tabular}
  \right.
\end{equation}
Note that $\bfv_i \in \bfH^{\dive}(\Omega_{\mcT_{h,i}})\cap \prod_{T\in\mcT_{h,i}} \bfH^{1}(T)$ with $\dive \bfv_i = -q_{h,i}$ in $\Omdhi$ and we have from the commuting diagram the identity (see \cite[Theorem 3.2 and Lemma 3.7]{Gatica14}):
\begin{align}
\dive\bfv_{h,i}= -q_{h,i} \quad \textrm{in $\Omdhi$.}
\label{eq: commuting property}
\end{align}
Thus, for $q_{h}=(q_{h,1},q_{h,2})$ and $\bfv_{h}=(\bfv_{h,1},\bfv_{h,2})$ we have
\begin{align*}
	\Bf(\bfv_{h},q_{h}) &= -(\dive \bfv_{h},q_{h})_{\Omega_1\cup\Omega_2} - s_b(\bfv_{h},q_{h})\\
	&= \|q_{h} \|_{\Omega_1\cup\Omega_2}^2 + s_b(q_h, q_h)
        \gtrsim \| q_{h} \|^2_h,
\end{align*}
with the last inequality being an application of Lemma \ref{lem:sp_ineq}. We are done if we can show that $\tn \bfv_h \tn_h \lesssim \| q_h\|_h$.

Recall the definition of the norm $\tn \bfv_{h} \tn_h^2$,
\begin{align*}
	\tn \bfv_{h} \tn_h^2 &= \underbrace{\|\permi^{1/2}\bfv_h\|^2_{\Omega_1\cup\Omega_2}+ \|\dive\bfv_h\|^2_{\Omdh} }_{\mathbf{I}}
	+
	\underbrace{s_{\bfu}(\bfv_h,\bfv_h)}_{\mathbf{II}}
	+
	\underbrace{ \|\permig^{1/2}\{\bfv_h\cdot\pmb{n}\}\|^2_{\interf} + \xi \|\permig^{1/2}\jump{\bfv_h\cdot\pmb{n}}\|^2_{\interf}}_{\mathbf{III}}.
\end{align*}

\textbf{Term I}. Recall our choice of $\bfv_{h,i}$ (see equation \eqref{eq:defvhi}). First note that using local interpolation estimates, Lemma \ref{lem:loc_interp} with $l=m=0$, we have 
\begin{align}\label{eq:intponv}
	\|\bfv_{h,i}\|^2_{T} &\lesssim \|\bfv_i-\bfv_{h,i}\|^2_{T}  + \|\bfv_i\|^2_{T} \lesssim h^2\|\nabla\bfv_i\|^2_{T} + \|\bfv_i\|^2_{T} \lesssim \|\bfv_i\|^2_{\bfH^1(T)}. 
\end{align}
Next using that $\bigcup_{T \in \mathcal{G}_h } T \subset \intfneig$ together with equations \eqref{eq: vpnorm} and \eqref{eq:vinorm} and that $q_{h,i}^e=0$ outside $\Omdhi$  we obtain
\begin{align}
\sum_{i=1}^2 \sum_{T\in\mcT_{h,i}} \|\bfv_i\|_{\bfH^1(T)} 
&
\lesssim
\|\bfv_p\|_{\bfH^1(\Omega)}+\sum_{i=1}^2\|\bfv_i^e\|_{\bfH^1( \intfneig)}
\lesssim \sum_{i=1}^2 \left(\|q_{h,i}\|_{\Omega_i}+\|q_{h,i}^e\|_{\intfneig} \right) \lesssim
\sum_{i=1}^2\|q_{h,i}\|_{\Omdhi}.
\label{eq:estimatvinq}
\end{align}
Applying \eqref{eq: commuting property}-\eqref{eq:estimatvinq} we get 
\begin{align}
\|\permi^{1/2}\bfv_h\|^2_{\Omega_1\cup\Omega_2}+ \|\dive\bfv_h\|^2_{\Omdh}  \lesssim \sum_{i=1}^2\sum_{T\in\mcT_{h,i}} \left(\|\bfv_{h,i}\|_{T}^2+\|\dive \bfv_{h,i}\|_{T}^2 \right) \lesssim  \|q_h\|^2_h.
\end{align}

\textbf{Term II}. Using a standard elementwise trace inequality followed by the standard inverse inequality to remove $j-1$ derivatives,  see e.g. \cite{BreSco}, we obtain 
\begin{align}\label{eq:stabbound}
& s_{\bfu}(\bfv_{h},\bfv_{h})=
\sum_{i=1}^2 \sum_{F \in \mathcal{F}_{h,i}} \tau_{u} \left( h \|\jump{\pmb{v}_{h,i}}\|^2_{F}  + \sum_{j=1}^{\kp+1} h^{2j + 1} \|\jump{D^j_{\bfn_F} \pmb{v}_{h,i}}\|^2_{F} \right)  \nonumber \\
& \quad  \lesssim  
\sum_{i=1}^2\sum_{T\in\mcT_{h,i}} \left( \| \bfv_{h,i}\|_{T}^2 + \sum_{j=1}^{\kp+1} h^{2j} h^{2-2j} \| \bfv_{h,i} \|_{\bfH^1(T)}^2 \right).
\end{align}
Local interpolation estimates, see Lemma \ref{lem:loc_interp}, yield
\begin{align}\label{eq:H1bound}
	\|\bfv_{h,i}\|^2_{\bfH^1(T)} &\lesssim \|\bfv_i-\bfv_{h,i}\|_{\bfH^1(T)}^2 + \|\bfv_i\|_{\bfH^1(T)}^2 \lesssim \|\bfv_i\|^2_{\bfH^1(T)}.
\end{align}
Applying this estimate together with \eqref{eq:estimatvinq} we get
\begin{align}
& s_{\bfu}(\bfv_{h},\bfv_{h}) \lesssim  
 \sum_{i=1}^2 \sum_{T\in\mcT_{h,i}} \| \bfv_i\|_{\bfH^1(T)}^2 
  \lesssim \|q_h\|^2_h.
\end{align}

\textbf{Term III}.  
We handle the interface terms by applying an elementwise trace inequality followed by the standard inverse inequality, the estimate in \eqref{eq:intponv}, using that $\cup_{T\in\mcGh} T \subset U_\varepsilon(\Gamma)$ and on this tubular neighborhood we apply the Poincaré inequality, Lemma \ref{lem:Poincare}. We finally use that $\varepsilon=ch$. Thus,   
\begin{align} 
  \|\bfv_{h,i}\cdot \pmb{n}\|^2_{\interf} &=\sum_{T\in\mcGh} \|\bfv_{h,i}\cdot \pmb{n}\|^2_{\interf\cap T}
  \lesssim h^{-1} \sum_{T\in\mcGh}\|\bfv_{h,i}\|^2_{T} 
        \lesssim h^{-1} \sum_{T\in\mcGh}(h^2\|\nabla\bfv^e_i\|^2_{T}+\|\bfv^e_i\|^2_{T}) \\ \nonumber
       	&\lesssim h\|\nabla\bfv_i^e\|^2_{U_\varepsilon(\Gamma)}+ h^{-1}\|\bfv_i^e\|^2_{U_\varepsilon(\Gamma)} \\ \nonumber
	&\lesssim \max(h,h^{-1}\varepsilon^2) \|\bfv_i^e\|^2_{\bfH^1(U_\varepsilon(\Gamma))} \lesssim \|q_{h,i}\|^2_{\Omdhi},
\end{align}
with the last inequality being due to \eqref{eq:vinorm} and since $q_{h,i}^e=0$ outside $\Omdhi$.

After splitting up averages and jump terms and applying the above estimate we get,
\begin{align}\label{eq:termI-IV}
	\|\permig^{1/2}\{\bfv_h\cdot\pmb{n}\}\|^2_{\interf} + \xi \|\jump{\permig^{1/2}\bfv_h\cdot\pmb{n}}\|^2_{\interf}
	 \lesssim \sum_{i=1}^2\|q_{h,i}\|^2_{\Omdhi} = \|q_h\|^2_h.
\end{align}
 
Combining the estimates for \textbf{Term I-III} yields the desired inequality,  
\begin{align*}
	\tn \bfv_{h} \tn_h^2 &\lesssim  \| q_{h} \|_h^2. 
\end{align*}    
\end{proof}

\subsection{Interpolation error estimates}\label{sec:interpol}
To define the interpolant we need extensions of functions in $\Omega_i$ to $\Omdhi$.
We need extension operators
\begin{equation}
	\bfE_i: \Hdivinto(\Omega_i)\to \Hdivinto(\RR^d),\
        \hat{E}_i: H^{\intorp}(\Omega_i)\to H^{\intorp}(\RR^d), 
\end{equation}
with $\intorp\geq 1$, that for $\bfv_i\in\Hdivinto(\Omega_i)$ and $q_i\in H^{\intorp}(\Omega_i)$ satisfy
\begin{align}
\|\bfE_i\bfv_i\|_{\Hdivinto(\RR^d)} &\lesssim \|\bfv_i\|_{\Hdivinto(\Omega_i)},
\quad (\bfE_i \bfv_i)|_{\Omega_i}=\bfv_i \text{ a.e.},
\label{eq:sobo_velbound}\\
\|\hat{E}_i  q_i\|_{H^{\intorp}(\RR^d)} &\lesssim \|q_i\|_{H^{\intorp}(\Omega_i)}, \quad
 (\hat{E}_i q_i)|_{\Omega_i}=q_i \text{ a.e.}, \label{eq:sobo_presbound}\\
\dive\circ \bfE_i &= \hat{E}_i \circ \dive.\label{eq:commdiagE}
\end{align}
The Sobolev-Stein extension operators satisfy these properties, see \cite[(3.16) and Corollary 4.1]{Hiptmair2012Universal}. 
We can then define

\begin{equation}
	\bfE: \bigoplus_{i=1}^2 \Hdivinto(\Omega_i) \ni (\bfv_1, \bfv_2) \to (\bfE_1 \bfv_1, \bfE_2 \bfv_2) \in \bigoplus_{i=1}^2 \Hdivinto(\RR^d),
\end{equation}

\begin{equation}
	\hat{E}: \bigoplus_{i=1}^2 H^{\intorp}(\Omega_i) \ni (q_1, q_2) \to (\hat{E}_1  q_1, \hat{E}_2  q_2) \in \bigoplus_{i=1}^2 H^{\intorp}(\RR^d),
\end{equation}
and finally the following interpolation and projection operators

\begin{equation}\label{eq:defintpv}
\interpv : \bigoplus_{i=1}^2 \Hdivinto (\Omega_i) \ni (\bfv_1, \bfv_2) \to \left(\pi_{h,1}^k \left(\bfE_1\bfv_1|_{\Omega_{\mcT_{h,1}}}\right), \pi_{h,2}^k \left(\bfE_2\bfv_2|_{\Omega_{\mcT_{h,2}}}\right) \right) \in \Vk, 
\end{equation}

\begin{equation}\label{eq:defintpp}
\interpp : \bigoplus_{i=1}^2 H^{\intorp}(\Omega_i) \ni (q_1, q_2) \to \left(\projpii{1} \left(\hat{E}_1q_1|_{\Omega_{\mcT_{h,1}}}\right), \projpii{2} \left(\hat{E}_2 q_2|_{\Omega_{\mcT_{h,2}}}\right) \right)\in \Qk. 
\end{equation}
\begin{lemma}\textbf{(Interpolation/projection estimates)}\label{lem:interp}
	Let $\bfu\in \bigoplus_{i=1}^2 \Hdivinto(\Omega_i)$ and $p\in \bigoplus_{i=1}^2 H^{\intorp}(\Omega_i)$, then
        \begin{align}
	 \tn \bfE \bfu-\interpv \bfu \tn_h &\lesssim h^{k+1/2}\sum_{i=1}^2\left(\|\bfu_i\|_{\bfH^{k+1}(\Omega_i)} + \|\dive\bfu_i\|_{H^{k+1}(\Omega_i)} \right),\label{eq:interpresu} \\
        \| \hat{E}p-\interpp p \|_h &\lesssim h^{k+1}\sum_{i=1}^2 \|p_i\|_{H^{k+1}(\Omega_i)}.\label{eq:interpresp}
        \end{align}
\end{lemma}
\begin{proof}
Let $\interr= \bfE \bfu-\interpv \bfu$. Recall that 
	\begin{align*}
		\tn \interr \tn_h^2 &= 
		\|\permi^{1/2}\interr \|^2_{\Omega_1\cup\Omega_2}  + \|\dive\interr \|^2_{\Omdh} +
		s_{\bfu}(\interr, \interr)	+
		\|\permig^{1/2}\{\interr \cdot\pmb{n}\}\|^2_{\interf} + \xi \|\permig^{1/2}\jump{\interr \cdot\pmb{n}}\|^2_{\interf}.
	\end{align*}
Using local interpolation estimates from Lemma \ref{lem:loc_interp} and the properties of the extension operators, \eqref{eq:sobo_velbound}-\eqref{eq:commdiagE}, we have
\begin{align}
\|\permi^{1/2}\interr_i\|^2_{\Omega_i}& \lesssim \sum_{T\in \mcTh}\|\permi^{1/2}\interr_i\|^2_{T} \lesssim h^{2(k+1)} \sum_{T\in \mcTh} |\bfE_i \bfu_i |^2_{k+1,T} \lesssim h^{2(k+1)} \| \bfu_i \|^2_{\bfH^{k+1}(\Omega_i)},\label{eq:interpboundI} \\
	\|\dive \interr_i \|^2_{\Omdhi} &= \sum_{T\in \mcTh} \|\dive \interr_i \|^2_{T} \lesssim h^{2(k+1)} \sum_{T\in \mcTh} |\dive \bfE_i\bfu_i |^2_{k+1,T} 
		=h^{2(k+1)} \sum_{T\in \mcTh} |\hat{E}_i \dive \bfu_i|_{k+1,T}^2  \nonumber \\
                &\lesssim h^{2(k+1)} \| \dive\bfu_i \|^2_{H^{k+1}(\Omega_i)}. \label{eq:interpboundI_2}
\end{align}

Next we look at the stabilisation term. We use an elementwise trace inequality, followed by local interpolation estimates ($\|D^j\interr_i\|_{T}\lesssim h^{k+1-j}|\bfE_i\bfu_i|_{k+1,T}$ from Lemma \ref{lem:loc_interp}), and the stability of the extension operator to get
\begin{align}
	s_{\bfu}(\interr, \interr)&=\sum_{i=1}^2 \sum_{F \in \mathcal{F}_{h,i}} \tau_{u} \left(\sum_{j=0}^{k+1} h^{2j + 1} \|\jump{D^j_{\bfn_F} \interr_i}\|^2_{F} \right) \lesssim \sum_{i=1}^2 \sum_{T \in \mcGh} \left( 
	\sum_{j=0}^{k+1} h^{2j+1}\left(h^{-1}\|D^j \interr_i\|_{T}^2+h\|D^{j+1} \interr_i\|_{T}^2 \right) \right) \nonumber \\
	&\lesssim h^{2(k+1)}\sum_{i=1}^2 \sum_{T\in \mcTh} | \bfE_i \bfu_i |^2_{k+1,T} \lesssim h^{2(k+1)}\sum_{i=1}^2 \| \bfu_i \|^2_{\bfH^{k+1}(\Omega_i)}. \label{eq:interpboundII}
\end{align}
For the interface terms we again use an elementwise trace inequality, followed by local interpolation estimates, and the stability of the extension operator. We get
 
\begin{align}
	\|\interr_i \cdot \pmb{n}\|^2_{\interf} &
	\lesssim \sum_{T\in\mcGh} \left(h^{-1} \|\interr_i\|^2_{T}
+h \|D \interr_i\|^2_{T} \right)
\nonumber \\
	&\lesssim h^{2k+1}\sum_{T\in\mcTh}|\bfE_i\bfu_i|^2_{k+1,T} 
	\lesssim h^{2k+1}\|\bfu_i\|_{\bfH^{k+1}(\Omega_i)}^2.
\end{align} 
Summing all terms together we get the desired estimate \eqref{eq:interpresu}. Since $\|\cdot \|_h$ is just the broken $L^2$-norm on the active meshes, the estimate \eqref{eq:interpresp} follows from standard estimates and using the stability of the extension operator $\hat{E}_i$, see \eqref{eq:sobo_presbound}. 
\end{proof}

\begin{remark}\textbf{The interpolation estimate}
	The interface terms are the only terms which contribute to the $h^{k+1/2}$ scaling of the interpolation estimate. For all other terms we get $h^{k+1}$. 
	We note that it is not possible to take the interpolation estimates of Lemma \ref{lem:loc_interp} beyond $|\bfu|_{k+1,T}$ in the hopes of picking up an extra power of $h$. For example for $k=0$ one can construct a linear function $\bfv$ on an element $T$ so that $\bfv \in \bfQ_1(T)\setminus\RT_0(T)$. If the local interpolation estimates of Lemma \ref{lem:loc_interp} could go to $|\cdot|_{2,T},$ then
	\begin{align}
		|\bfv-\interpv \bfv|_{0,T}\lesssim h^2|\bfv|_{2,T}=0,
	\end{align}
	which is a contradiction due to the construction of $\bfv$ and that $\interpv \bfv$ is a function in $\RT_0(T)$.
\end{remark} 

\subsection{Consistency and continuity}\label{sec:consistency}
Before we are ready to show a priori error estimates we need to show a result on consistency of Method 2 and that the bilinear forms $\Af$ and $\Bf$ are continuous.

\begin{lemma}(\textbf{Consistency})\label{lem:consist}
Let $\bfu \in \bfV \cap \bigoplus_{i=1}^2 \Hdivinto(\Omega_i)$ and $p \in \bigoplus_{i=1}^2 H^{k+1}(\Omega_i)$ be the solution to \eqref{eq:weakdarcy}-\eqref{eq:weakdarcy2}. Let $(\bfu_h,p_h)\in\Vk\times\Qk$ be the solution to \eqref{eq:discretedarcy2}-\eqref{eq:discretedarcy2b}. Then,
\begin{align}
	\Af(\bfu_h-\bfE \bfu,\pmb{v}_h) + \Bf(\bfv_h, p_h-\hat{E} p) &=-s_{\bfu}(\bfE \bfu,\bfv_h),
        \quad \forall \bfv_h \in \Vk,\label{eq:consisteq1}\\
	\Bf(\bfu_h-\bfE \bfu,q_h) &= 0, \quad \forall q_h \in \Qk,
\end{align}
with
\begin{equation}\label{eq:suforEu}
  s_{\bfu}(\bfE \bfu,\bfv_h)= \sum_{i=1}^2 \sum_{F \in \mathcal{F}_{h,i}}
\tau_{\bfu} h^{2k+3}  (\jump{D^{k+1}_{\bfn_F} \bfE_i \bfu_i }, \jump{D^{k+1}_{\bfn_F} \bfv_{h,i}})_{F}.
\end{equation}

\end{lemma}
 
\begin{proof} 
  Note that since $\hat{E_i}p_i \in H^{k+1}(\RR^d)$ and $\bfE_i \bfu_i \in \Hdivinto(\RR^d)$ we have $\jump{D^{j} \hat{E_i}p_i}=0$, $\jump{D^j (\dive (\bfE_i \bfu_i)}=0$, and $\jump{D^{j}_{\bfn_F} (\bfE_i \bfu_i)^l}=0$   for $0\leq j \leq k$, ($l=1, \cdots, d$). Hence $s_b(\bfv_h,\hat{E}p)=s_b(\bfE \bfu, q_h)=0$ and consequently $\Bf(\bfE \bfu, q_h)=b(\bfu, q_h)$ and $\Bf(\bfv_h, \hat{E} p)=b(\bfv_h, p)$. Furthermore, equation \eqref{eq:suforEu} holds. Following the derivation of the weak form \eqref{eq:weakdarcy}-\eqref{eq:weakdarcy2}, we get 
\begin{align}
  \Af(\bfu_h-\bfE \bfu,\bfv_h) + \Bf(\bfv_h, p_h-\hat{E} p) &= \FF (\pmb{v}_h)-\left(a(\bfu, \bfv_h)+\sum_{i=1}^2 \sum_{F \in \mathcal{F}_{h,i}}
  \tau_{\bfu} h^{2k+3} ([D^{k+1}_{\bfn_F} \bfE_i \bfu_i ], [D^{k+1}_{\bfn_F} \bfv_{h,i}])_{F}  + b(\bfv_h, p) \right), \nonumber \\
  	\Bf(\bfu_h-\bfE \bfu,q_h) &= 0, \nonumber
\end{align}
and the result follows. 
\end{proof}

\begin{lemma}\textbf{(Continuity)}\label{lem:continuity}
The bilinear forms $\Af$ and $\Bf$ are continuous
	\begin{align}
		\Af(\bfv,\bfw) &\lesssim \tn \bfv \tn_h \tn \bfw \tn_h, \quad \forall \bfv, \bfw \in \left(\bigoplus_{i=1}^2\Hdivinto(\Omdhi)+\Vk \right),  \\ 
		\Bf(\bfv,q) &\lesssim \tn \bfv \tn_h\|q\|_h, \quad \forall \bfv \in \left(\bigoplus_{i=1}^2\Hdivinto(\Omdhi)+\Vk \right), \ q \in \bigoplus_{i=1}^2\left(H^{k+1}(\Omdhi) + \Qki \right).
	\end{align}
\end{lemma}
\begin{proof}
	Applying Cauchy-Schwartz we see that $\Af(\bfv,\bfw) \lesssim \Af(\bfv,\bfv)^{1/2}\Af(\bfw,\bfw)^{1/2}$ whereby the first inequality follows. To show that $\Bf$ is continous we again apply Cauchy-Schwartz,
	\begin{align}\label{eq:contB}
		\Bf(\bfv,q) &= -(\dive\bfv,q)_{\Omega_1\cup\Omega_2}-s_b(\bfv,q) \lesssim \| \dive\bfv \|_{\Omega_1\cup\Omega_2} \|q\|_{\Omega_1\cup\Omega_2} + |s_b(\bfv,q)| \nonumber \\
		&\lesssim \| \dive\bfv \|_{\Omega_1\cup\Omega_2} \|q\|_{\Omega_1\cup\Omega_2}  + s_p(\dive\bfv,\dive\bfv)^{1/2}s_p(q,q)^{1/2} \nonumber \\ 
		&\lesssim\| \dive\bfv \|_{h}\|q\|_{h} \lesssim \tn\bfv\tn_h\|q\|_h.
	\end{align}
        Note that for the stabilization term we have used that for $q_h \in \Qk$ an elementwise trace inequality followed by the standard inverse inequality to remove $j$ derivatives yields
\begin{align}\label{eq:stabbound}
& s_{p}(q_{h},q_{h})=
\sum_{i=1}^2 \sum_{F \in \mathcal{F}_{h,i}} \tau_p \left( h \|\jump{q_{h,i}}\|^2_{F}  + \sum_{j=1}^{k} h^{2j + 1} \|\jump{D^j q_{h,i}}\|^2_{F} \right)  \nonumber \\
& \quad  \lesssim  
\sum_{i=1}^2 \left( \| q_{h,i}\|_{\Omdhi}^2 + \sum_{T\in\mcT_{h,i}}\sum_{j=1}^{k} h^{2j} h^{-2j} \| q_{h,i} \|_{T}^2 \right)\lesssim \| q_{h} \|_{h}^2.
\end{align}
Furthermore, for $q_i \in H^{k+1}(\Omdhi)$ and $\dive \bfv_i \in  H^{k+1}(\Omdhi)$, we have $\jump{D^j q_{i}}=0$ and $\jump{D^j (\dive \bfv_{i})}=0$ for $0\leq j \leq k$. Thereof, $s_p(q,q)^{1/2} \lesssim \|q\|_h$ and $s_p(\dive\bfv,\dive\bfv)^{1/2} \lesssim \| \dive\bfv \|_{h}$ in equation \eqref{eq:contB}.    

\end{proof}

\subsection{A priori error estimates}\label{sec:apriori}
When the exact divergence is in the pressure space, $\dive\bfE\bfu\in\Qk$, the saddle-point structure of our problem allows us to mimic the a priori proof for standard FEM, see \cite[Theorem 2.6]{Gatica14}. However, for the general case 
we have to couple the pressure error with the velocity error. We introduce for $\bfu_h,\bfv_h\in\Vk,p_h,q_h\in\Qk,$
\begin{align}
	\Cf((\bfu_h,p_h),(\bfv_h,q_h)) &:= \Af(\bfu_h,\bfv_h) +\Bf(\bfv_h,p_h)+
	\Bf(\bfu_h,q_h) = \FF (\pmb{v}_h) + \GG(q_h), \\
	\tn (\bfv_h,q_h) \tn^2_h &:= \tn \bfv_h\tn^2_h + \|q_h\|^2_h,
\end{align}
next we prove a stability result, and finally we show optimal a priori error estimates. 
\begin{lemma}\label{lem:big_infsup}
	For any $(\bfrho_h,\sigma_h)\in \Vk\times\Qk$ there exists $(\bfw_h,r_h)\in \Vk\times\Qk$ such that 
	\begin{align}
          \Cf((\bfrho_h,\sigma_h),(\bfw_h,r_h))\gtrsim \tn (\bfrho_h,\sigma_h) \tn_h^2, \quad
          \tn (\bfw_h,r_h) \tn_h \lesssim \tn (\bfrho_h,\sigma_h) \tn_h.
	\end{align}
\end{lemma}  
\begin{proof}
	Let $-\bfv_h\in\Vk$ be the element attaining the inf-sup condition in Proposition \ref{prop:infsup}, that is $\Bf(-\bfv_h,\sigma_h)\gtrsim \|\sigma_h\|^2_h$, and $\tn\bfv_h\tn_h \lesssim \|\sigma_h\|_h$. Using also continuity of the bilinear form $\Af$ and Young's inequality we have
	\begin{align}
		\Cf((\bfrho_h,\sigma_h),(-\bfv_h,0)) &= -\Af(\bfrho_h,\bfv_h)+\Bf(-\bfv_h,\sigma_h) \nonumber \\
		&\gtrsim -\tn\bfrho_h\tn_h\tn\bfv_h\tn_h + \|\sigma_h\|^2_h \gtrsim -\tn\bfrho_h\tn_h\|\sigma_h\|_h + \|\sigma_h\|^2_h \nonumber \\
		&\gtrsim -\epsilon^{-1}\tn\bfrho_h\tn^2_h + (1-\epsilon)\|\sigma_h\|^2_h.
	\end{align}
Using Lemma \ref{lem:sp_ineq} we have $\Bf(\bfrho_h, -\dive\bfrho_h)=\|\dive \bfrho_h\|^2_{\Omega_1 \cup \Omega_2}+ s_p(\dive \bfrho_h, \dive \bfrho_h) \gtrsim \|\dive \bfrho_h\|^2_{h}$ and hence, 
	\begin{align}
		\Cf((\bfrho_h,\sigma_h),(\bfrho_h,-\dive\bfrho_h-\sigma_h)) &= \Af(\bfrho_h,\bfrho_h)+\Bf(\bfrho_h,\sigma_h)+\Bf(\bfrho_h, -\dive\bfrho_h-\sigma_h)\nonumber \\
		&\gtrsim \tn\bfrho_h\tn^2_{\Af} +  \|\dive \bfrho_h\|^2_{h}= \tn\bfrho_h\tn^2_{h}.
	\end{align}
        Choosing $(\bfw_h,r_h)=(\bfrho_h-\delta\bfv_h,-\dive\bfrho_h-\sigma_h)$ with $\delta>0$ we have
	\begin{align}
		\Cf((\bfrho_h,\sigma_h),(\bfrho_h-\delta\bfv_h,-\dive\bfrho_h-\sigma_h))
		&\gtrsim (1-\delta\epsilon^{-1}) \tn\bfrho_h\tn^2_{h} + \delta(1-\epsilon)\|\sigma_h\|^2_h \gtrsim \tn (\bfrho_h,\sigma_h) \tn^2_h, 
	\end{align}
        to get the last inequality we can for instance choose $\epsilon=1/2$ and $\delta=1/3$. Using that $\| \dive\bfrho_h \|_h \lesssim \tn \bfrho_h \tn_h$ and $\tn \bfv_h \tn_h \lesssim \|\sigma_h\|_h$ we also have that
	\begin{align}
          \tn (\bfw_h,r_h) \tn_h \lesssim \tn \bfrho_h-\delta\bfv_h\tn_h+ \| \dive\bfrho_h+\sigma_h\|_h 
          \lesssim \tn (\bfrho_h,\sigma_h) \tn_h.
	\end{align}
\end{proof}

\begin{thm}(\textbf{A priori error estimate})
	Let $\bfu \in \bfV\cap \bigoplus_{i=1}^2 \Hdivinto(\Omega_i)$ and $p \in \bigoplus_{i=1}^2 H^{k+1}(\Omega_i)$ be the solution to \eqref{eq:weakdarcy}-\eqref{eq:weakdarcy2}. Let $(\bfu_h,p_h)\in\Vk\times\Qk$ be the solution to \eqref{eq:discretedarcy2}-\eqref{eq:discretedarcy2b}. Then,
	\begin{align}
		\tn \bfE \bfu-\bfu_h \tn_h
+	\| \hat{E}p-p_h \|_h 
&\lesssim h^{k+1/2}\sum_{i=1}^2\left(\|\bfu_i\|_{\bfH^{k+1}(\Omega_i)} + \|\dive\bfu_i\|_{H^{k+1}(\Omega_i)}+\|p_i\|_{H^{k+1}(\Omega_i)}\right).
\end{align}
\end{thm}
\begin{proof}
Adding and subtracting the interpolant and using the triangle inequality we have,
\begin{align}
\tn \bfE\bfu-\bfu_h \tn_h \leq  \tn \bfE\bfu-\interpv\bfu \tn_h + \tn \bfu_h-\interpv\bfu \tn_h, \\
\|\hat{E}p-p_h\|_h\leq \|\hat{E}p-\interpp p\|_h + \|p_h-\interpp p\|_h.  
\end{align}
We now want to estimate $\tn \bfu_h-\interpv\bfu \tn_h +\|p_h-\interpp p\|_h$ in terms of $\tn \bfE\bfu-\interpv\bfu \tn_h+ \|\hat{E}p-\interpp p\|_h$. Applying Lemma \ref{lem:big_infsup} to $(\bfrho_h, \sigma_h)=(\bfu_h-\interpv\bfu, p_h-\interpp p)$, we have that there exist $(\bfw_h,r_h)\in \Vk\times\Qk$ such that 
\begin{align}\label{eq:stabC}
  \tn (\bfrho_h,\sigma_h) \tn_h^2 \lesssim  \Cf((\bfrho_h,\sigma_h),(\bfw_h,r_h)), \quad 
\tn (\bfw_h,r_h) \tn_h \lesssim \tn (\bfrho_h,\sigma_h) \tn_h.
\end{align}
Consistency (Lemma \ref{lem:consist}) and continuity of the bilinear forms $\Af$ and $\Bf$ (Lemma~\ref{lem:continuity}) together with Cauchy-Schwartz inequality ($s_{\bfu}(\bfE \bfu,\bfw_h)\leq s_{\bfu}(\bfE \bfu,\bfE \bfu)^{1/2} s_{\bfu}(\bfw_h,\bfw_h)^{1/2}$) yield
\begin{align}
		\left|\Cf((\bfrho_h,\sigma_h),(\bfw_h,r_h))  \right|&
		= \left|\Af(\bfu_h-\interpv\bfu,\bfw_h) +\Bf(\bfw_h,p_h-\interpp p)+
		 \Bf(\bfu_h-\interpv \bfu,r_h)\right| \nonumber \\
	  	&= \left|\Af(\bfE\bfu-\interpv\bfu,\bfw_h) +\Bf(\bfw_h,\hat{E}p-\interpp p)+
	 	 \Bf(\bfE\bfu-\interpv \bfu,r_h) 
                 +s_{\bfu}(\bfE \bfu,\bfw_h) \right| \nonumber \\
 		        &\lesssim \left(\tn \bfE\bfu-\interpv\bfu \tn_h+\|\hat{E}p-\interpp p\|_h +  \left(\sum_{i=1}^2 \sum_{F \in \mathcal{F}_{h,i}} h^{2k+3} \| D^{k+1}_{\bfn_F} \interpv \bfu\|_F^2 \right)^{1/2} \right)\tn (\bfw_h,r_h) \tn_h.
\nonumber 
\end{align} 
In the last step we used that all terms in the stabilization $s_{\bfu}(\bfE \bfu,\bfE \bfu)$ vanish except the last term since $\jump{D^{j}_{\bfn_F} (\bfE_i \bfu_i)^l}=0$ for $0\leq j \leq k$, ($l=1, \cdots, d$) and we added and subtracted the interpolant. Using \eqref{eq:stabC} we end up with the estimate
\begin{align}
  \tn (\bfrho_h,\sigma_h) \tn_h^2  \lesssim  \left(\tn \bfE\bfu-\interpv\bfu \tn_h+\|\hat{E}p-\interpp p\|_h  +\left(\sum_{i=1}^2 \sum_{F \in \mathcal{F}_{h,i}} h^{2k+3} \| D^{k+1}_{\bfn_F} \interpv \bfu\|_F^2 \right)^{1/2} \right)\tn (\bfrho_h,\sigma_h) \tn_h.
\end{align} 
Using a standard elementwise trace inequality, an inverse estimate, local interpolation estimates (Lemma \ref{lem:loc_interp}), and the stability of the extension operator we have  
\begin{align} \label{eq:consterr}
   \sum_{F \in \mathcal{F}_{h,i}} h^{2k+3} \| D^{k+1}_{\bfn_F} \interpv \bfu\|_F^2 &\lesssim  \sum_{T \in \mathcal{T}_{h,i}} h^{2k+2} \| D^{k+1} \interpv \bfu\|_T^2   \nonumber \\
  &\lesssim
   \sum_{T \in \mathcal{T}_{h,i}} h^{2k+2} \left( \| D^{k+1} (\interpv \bfu-\bfE\bfu)\|_T^2+  \| D^{k+1} \bfE\bfu \|_T^2 \right)  \lesssim
  h^{2k+2} \|\bfu_i\|_{\bfH^{k+1}(\Omega_i)}^2.
\end{align}

Finally, using the interpolation and projection error estimates in Lemma \ref{lem:interp} and \eqref{eq:consterr} we get the result:
\begin{align*}
\tn \bfE\bfu-\bfu_h \tn_h+\|\hat{E}p-p_h\|_h&\lesssim  \tn \bfE\bfu-\interpv\bfu \tn_h +\|\hat{E}p-\interpp p\|_h +\left(\sum_{i=1}^2 \sum_{F \in \mathcal{F}_{h,i}} h^{2k+3} \| D^{k+1}_{\bfn_F} \interpv \bfu\|_F^2 \right)^{1/2}  \\
&\lesssim h^{k+1/2}\sum_{i=1}^2\left(\|\bfu_i\|_{\bfH^{k+1}(\Omega_i)} + \|\dive\bfu_i\|_{H^{k+1}(\Omega_i)}+\|p_i\|_{H^{k+1}(\Omega_i)}\right).
\end{align*}
\end{proof}

\section{Numerical Experiments}\label{sec:numex}
We consider two numerical examples which serve as benchmarks for the presented methods. Our results are based on the element pairs $\RT_0\times Q_0$, $\BDM_1\times Q_0$, and $\RT_1\times Q_1.$  We compare both the order of accuracy of the presented schemes and the scaling of the condition number of the associated linear system matrices with mesh size.
We study the convergence of the approximated velocity $\bfu_h$ and the pressure $p_h$ in the $L^2$-norm and the convergence for the divergence $\dive \bfu_h$ both in the $L^2$-norm and the $L^\infty$-norm. To compute the  $L^\infty$-norm of a function we take the maximum value of the function at a number of quadrature points. We use a high order quadrature rule in each element and on the boundary of each element.

The constants in the stabilization terms are all set to be equal to $1$, except for Example 2 with pair $\RT_0\times Q_0$ where they are set to $0.1$.  These constants are not optimally chosen and one can probably choose the constants so that the constant in the error or the condition number improves.

In general we do not have the exact interface $\Gamma$  but rather an approximation $\Gamma_h$. Consequently, we also have approximations $\Omega_{h,i}$ of $\Omega_i$, $i=1,2$.  Thus, in the implementation of all the schemes, we have replaced  $\Omega_i$ by $\Omega_{h,i}$ and $\Gamma$ by $\Gamma_h$.

The code used to run the numerical simulations is an in-house code written in C++ based on Freefem++ \cite{Freefem}. We use a direct solver to solve all linear systems. The source code is uploaded to this Github repository \cite{repo}. A Python notebook which runs our numerical examples can also be accessed. See \cite{Frachon2022Thesis} for some explanations of the data structures involved in the source code, along with explanations of how to use them.

\subsection{Example 1 ($\dive \bfu$ in the pressure space)} 
\begin{figure}[h!]
	\centering
	\begin{subfigure}[b]{0.45 \textwidth}  	
		\centering	
		\includegraphics[scale=0.15]{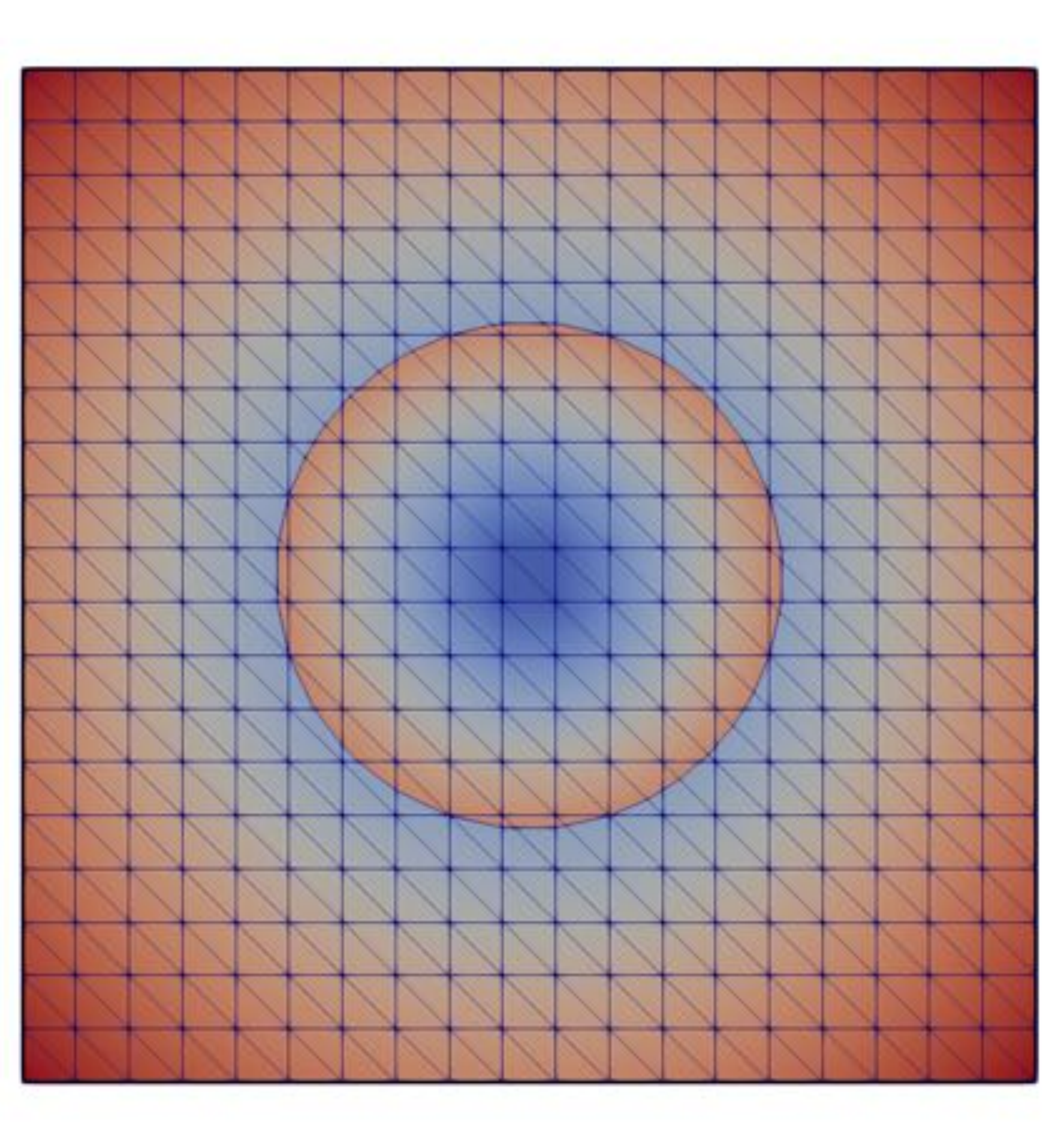}  
	\end{subfigure}
	\hfill
	\begin{subfigure}[b]{0.45\textwidth}  	
		\centering	
		\includegraphics[scale=0.2]{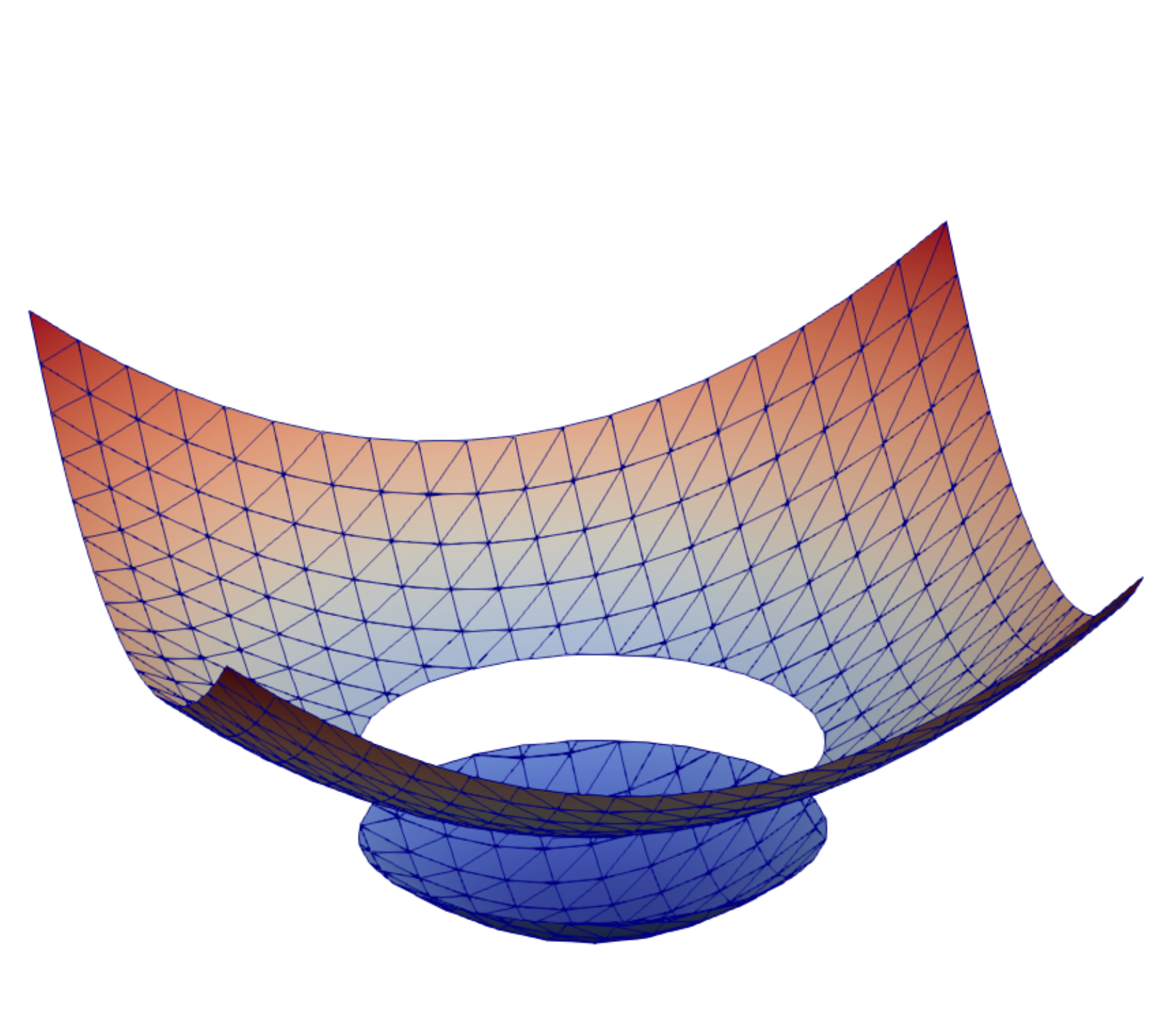}
	\end{subfigure} \\
		\begin{subfigure}[b]{0.45 \textwidth}  	
		\centering	
		\includegraphics[scale=0.2]{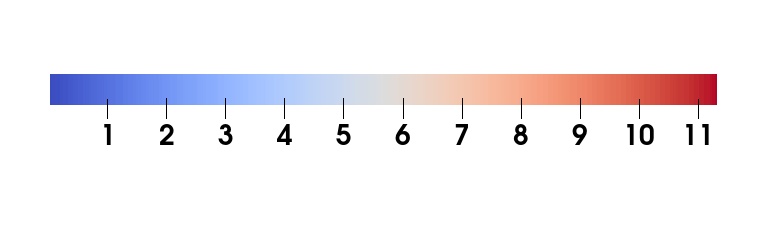}  
	\end{subfigure}
	\hfill
	\begin{subfigure}[b]{0.45\textwidth}  	
		\centering	
		\includegraphics[scale=0.24]{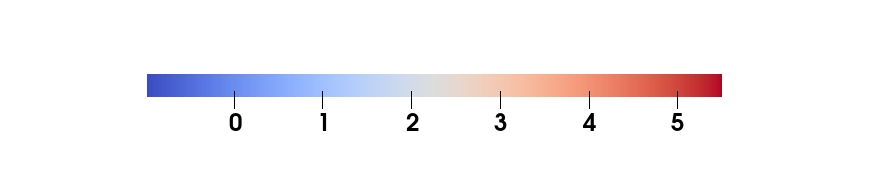}
	\end{subfigure} 
	\caption{   
	Example 1: Solution obtained using Method 2 with the element pair $\RT_1\times Q_1$ and a macro-element partitioning with $\delta = 0.25$. Left: Magnitude of the approximated velocity. Right: Approximated pressure. }
	\label{fig:example 1 solution}
\end{figure}

We consider the first example in \cite{DanSco12}. 
The computational domain is $\Omega = [0 , 1 ]\times  [0 , 1 ]$ and the interface is the circle defined as $$\Gamma = \left\{ (x,y) \in \RR^2 : r:=\sqrt{(x-1/2)^2+(y-1/2)^2} = R \right\} $$ with $R=0.25$. We set $\xi = 1/8,\ \permi=1$, $\permig = 2R/3$, the interface pressure $\hat{p}=19/12$, $\bff=0$, and 
\[ g = \begin{cases} -2/R^2 \text{ in } \Omega_1,\\-4/R^2 \text{ in } \Omega_2. \end{cases}\] 
We have natural boundary conditions, i.e. $\boundu=\emptyset$, with $p_B = r^2/(2R^2)+3/2$. 
The exact solution is therefore given by
\begin{align*}
    p=\begin{cases} r^2/(2R^2)+3/2 \text{ in } \Omega_1,\\ r^2/R^2 \text{ in } \Omega_2, \end{cases} & \text{ and }
    \pmb{u} = -\nabla p = \begin{cases} -1/R^2(x-1/2,y-1/2) \text{ in } \Omega_1,\\ -2/R^2(x-1/2,y-1/2) \text{ in } \Omega_2. \end{cases}
\end{align*}
We show the approximated solution obtained with Method 2 and the element pair $\RT_1\times Q_1$ in Figure \ref{fig:example 1 solution}. 

\begin{figure}[h]
	\centering 
	\includegraphics[scale=0.58]{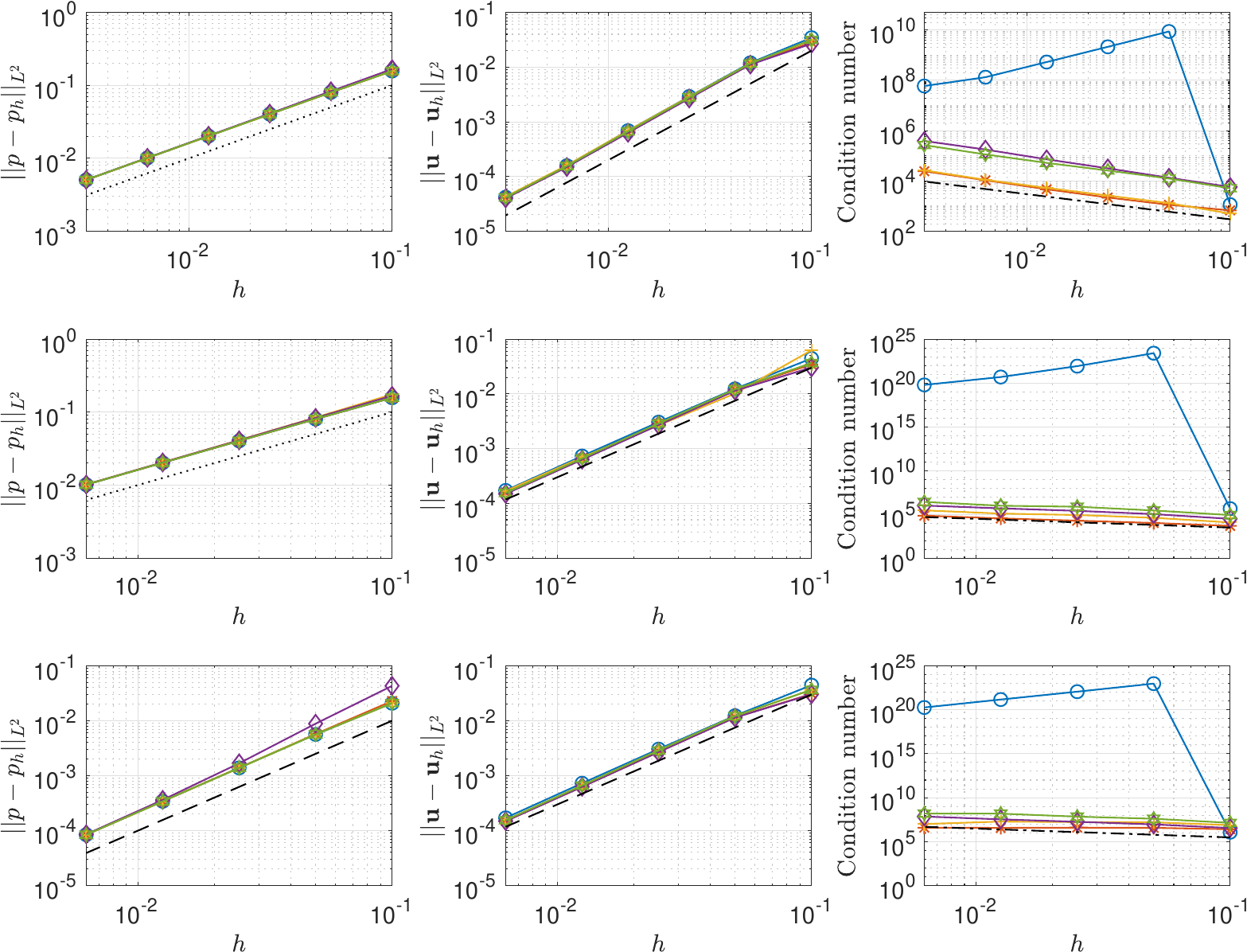} \\
	\includegraphics[scale=0.35]{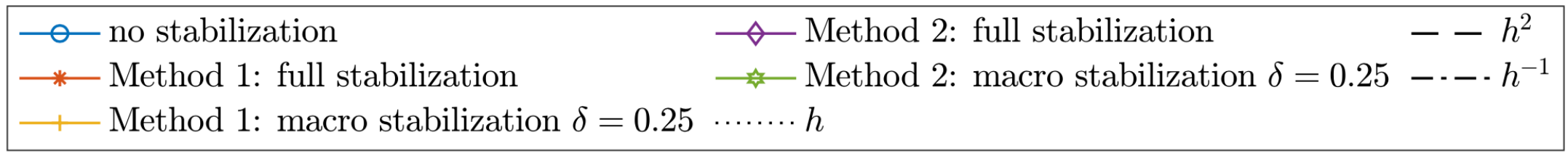}
	\caption{Example 1: Results using element pairs $\RT_0\times Q_0$ (first row), $\BDM_1\times Q_0$ (second row) and $\RT_1\times Q_1$ (third row). 
		Left: The $L^2$-error of the pressure versus mesh size h. Middle:  The $L^2$-error of the velocity field versus mesh size h. Right: The spectral condition number versus mesh size h. }     
	\label{fig:plot example1 convergence}  	
\end{figure}

We compare the two stabilized unfitted discretizations with the unstabilized discretization (which corresponds to taking the stabilization parameters equal to zero, i.e., $\tau_{\bfu}=\tau_{p}=0$ in Method 1), using the element pairs $\RT_0\times Q_0$, $\BDM_1\times Q_0$, and $\RT_1\times Q_1$. 
In Figure~\ref{fig:plot example1 convergence} we show the $L^2$-error of the approximated pressure and the velocity field versus mesh size h, and the spectral condition number of the linear system matrix associated with the methods presented in Section~\ref{sec:Nummeth}. 

\begin{figure}[h!]
	\centering
	\begin{subfigure}[b]{0.45 \textwidth}  	
		\centering	
		\includegraphics[scale=0.35]{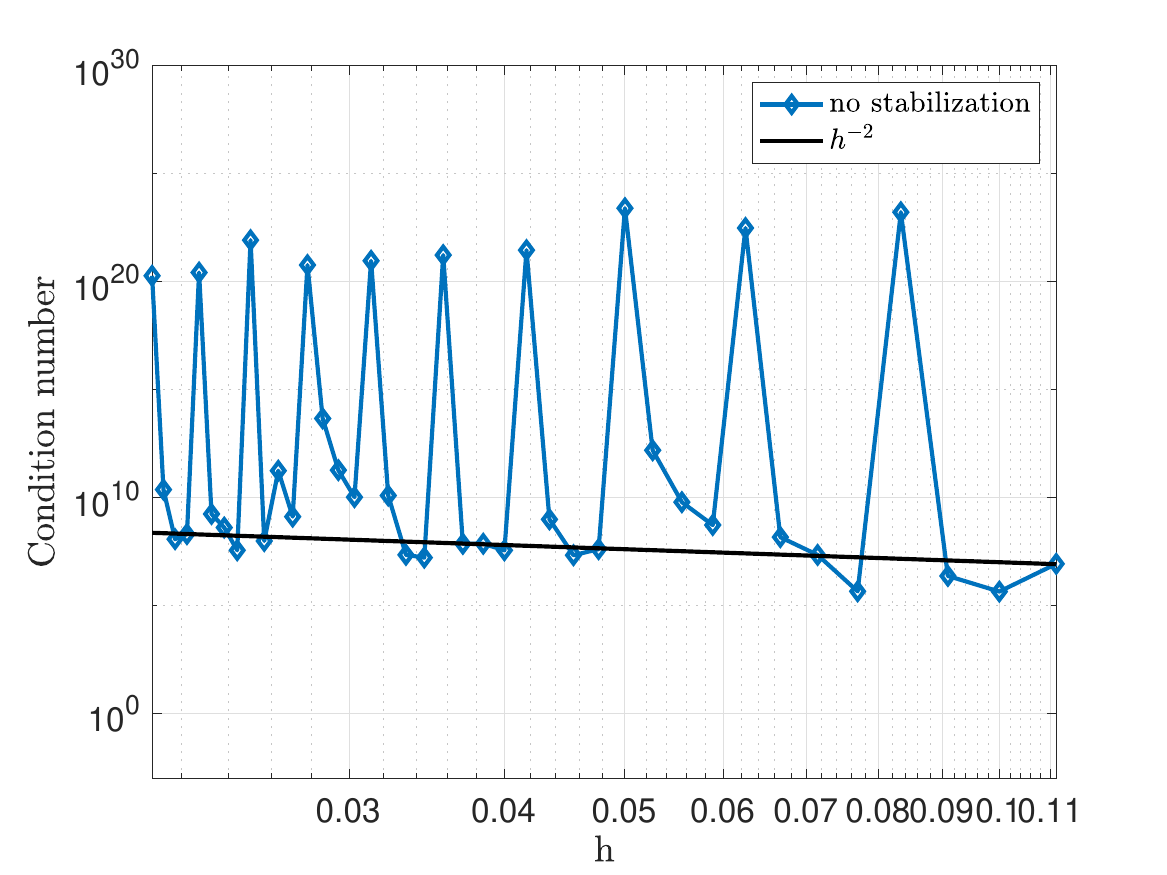}  
	\end{subfigure}
	\begin{subfigure}[b]{0.45 \textwidth}  	
		\centering	
		\includegraphics[scale=0.35]{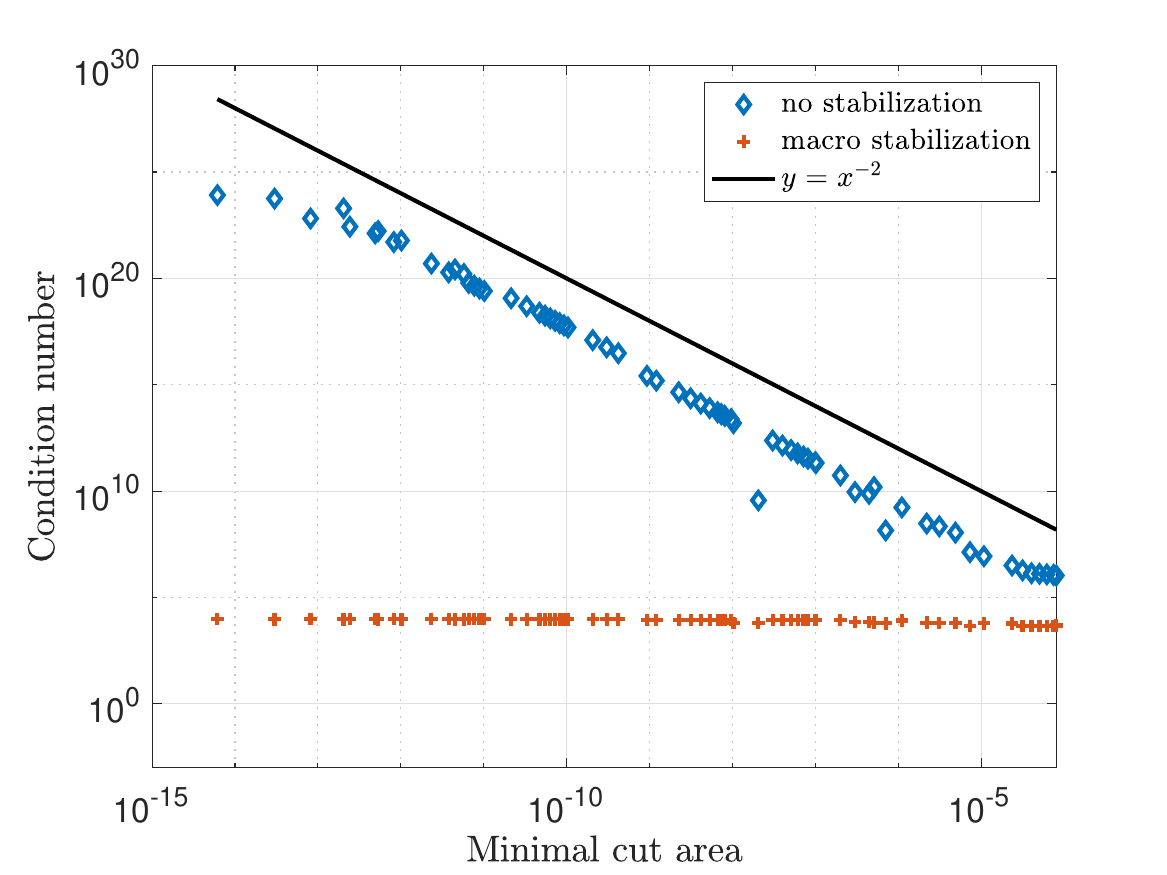}  
	\end{subfigure}
	\\
	\caption{Example 1:  Left: The spectral condition number for different mesh sizes.  Right: Condition number as a function of the smallest area of a cut region, $\min_{T \in \mcGh} |T\cap \Omega_i|$. Blue diamonds represent the unstabilized method and red stars represent Method 2 with macro stabilization.
	\label{fig:condition number step}}
\end{figure}

We obtain optimal convergence for the pressure with all the presented methods. For the velocity field we observe convergence of order $2$ for all methods and all elements. With the stabilized methods we obtain as expected optimal scaling of the condition number, that is $O(h^{-1})$. However, for the unstabilized method, the condition number is not controlled and may become arbitrarily large depending on the interface position. In Figure~\ref{fig:condition number step} we show how the condition number for the unstabilized method increases as the smallest area of the cut triangles, $\min_{T \in \mcGh} |T\cap \Omega_i|$,  decreases while the condition number stays constant for the stabilized method. We also show the condition number of the unstabilized method for many different mesh sizes in order to illustrate that the condition number varies depending on the interface position relative the mesh.

\begin{figure}[h]
\centering
    \begin{subfigure}[b]{0.45 \textwidth}  	
	\centering	
 	 \includegraphics[scale=0.17]{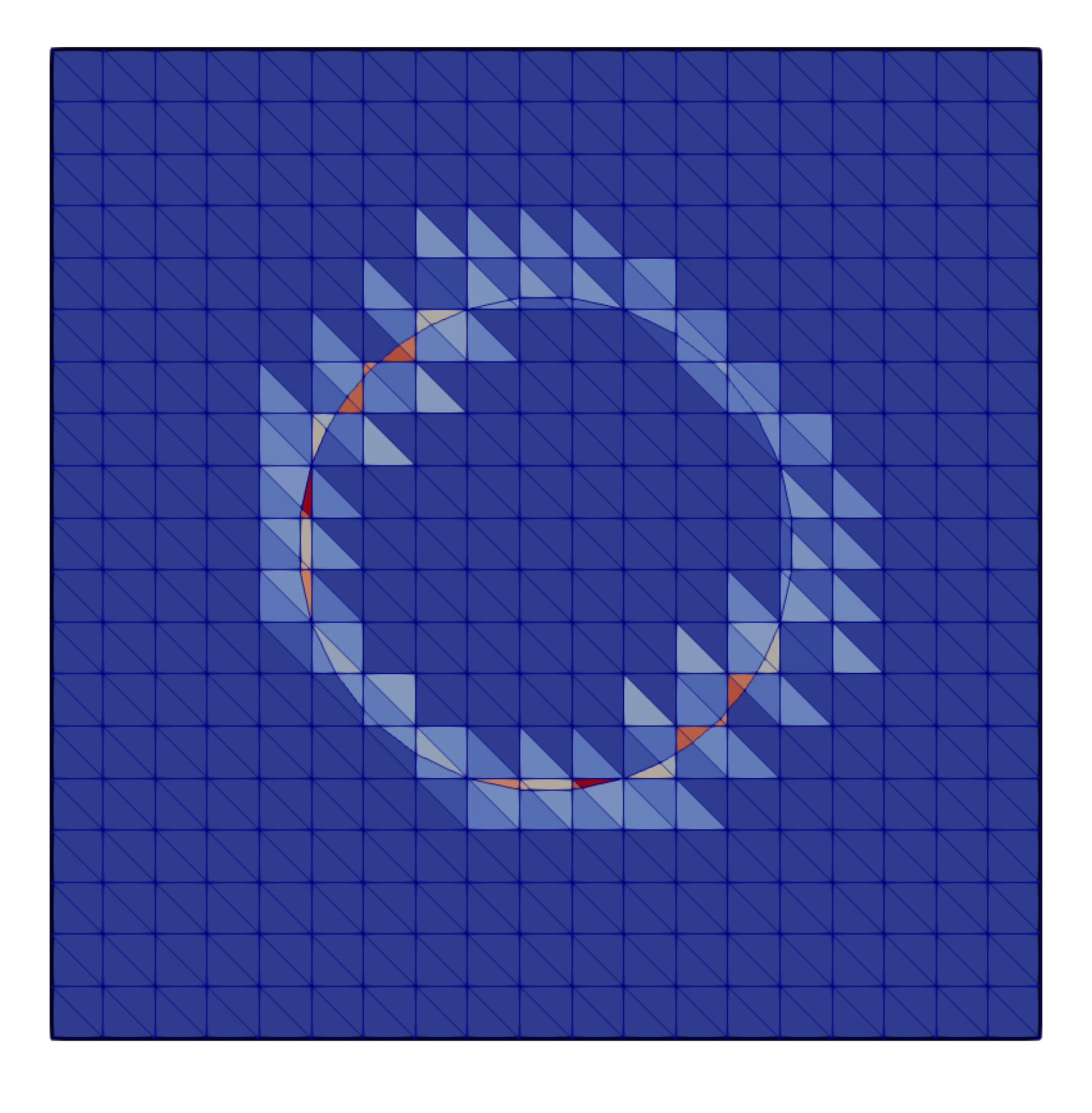}
    \end{subfigure} 
    \hfill
    \begin{subfigure}[b]{0.45\textwidth}  	
	\centering	
	\includegraphics[scale=0.17]{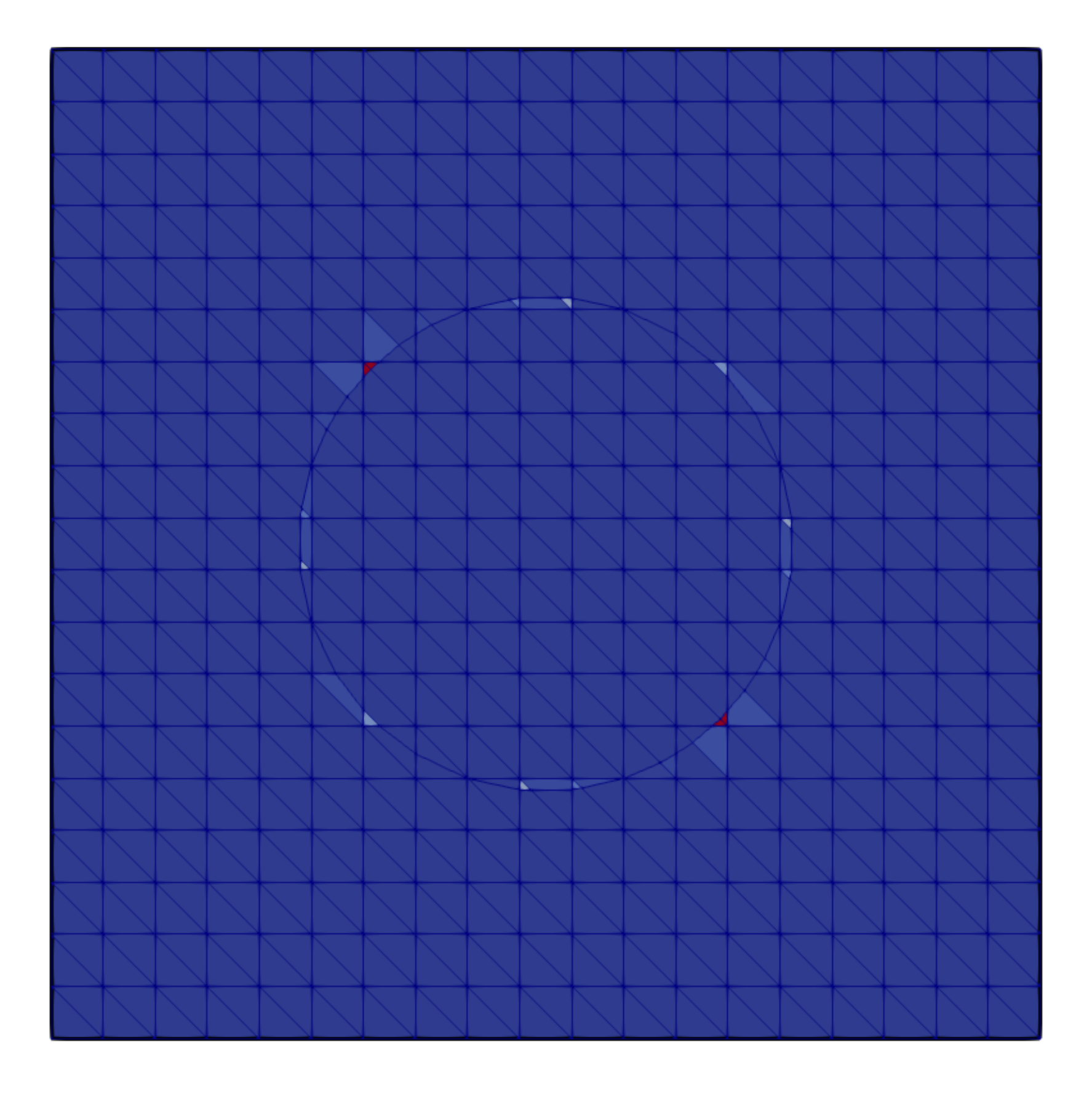}
    \end{subfigure}
\begin{subfigure}[b]{0.45\textwidth}  	
	\centering	
	\includegraphics[scale=0.25]{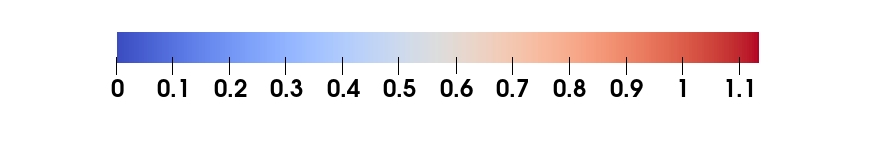}
    \end{subfigure}
         \caption{Example 1: The max-norm error $||\dive\pmb{u}-\dive\pmb{u}_h||_{L^{\infty}}$ using Method 1 with $\RT_0\times Q_0$. Left: The error using full stabilization. Right: The error using macro stabilization with  $\delta = 0.25$ where we get 6 small elements in $\Omega_1$ and 24 small elements in $\Omega_2$. 
           \label{fig:divergence error macro/full stabilization example 1}      
         }     	
\end{figure}

\begin{figure}[h]
    \centering
        \begin{subfigure}[b]{0.45 \textwidth}  	
	\centering	
    \includegraphics[scale=0.17]{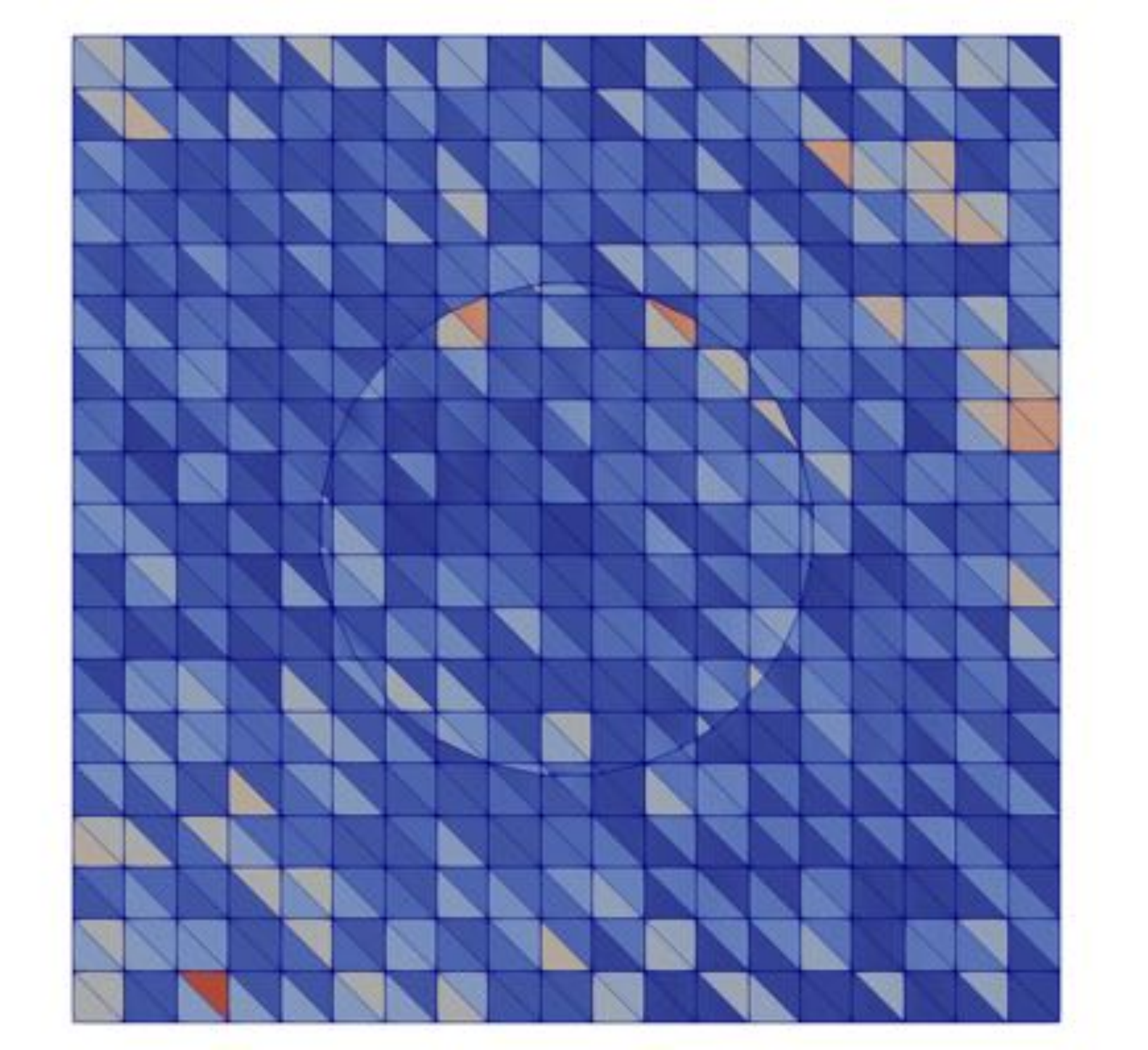}  
      \end{subfigure}\\
      \vspace{-0.5cm}
              \begin{subfigure}[b]{0.45 \textwidth}  	
	\centering	
    \includegraphics[scale=0.2]{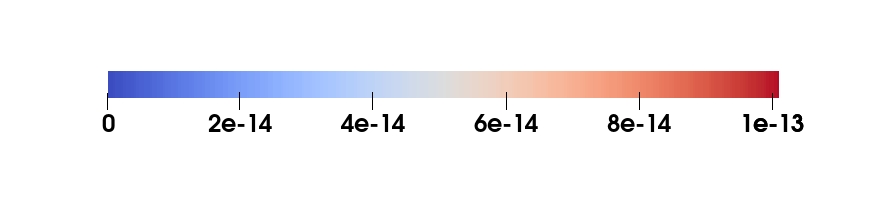}  
      \end{subfigure}
    \caption{Example 1: The max-norm error $||\dive\pmb{u}-\dive\pmb{u}_h||_{L^{\infty}}$ using Method 2 and full stabilization with $\RT_0\times Q_0$. }
    \label{fig:divergence error method2 example 1}
\end{figure}

We show the $L^{\infty}$-error of the divergence, $||\dive\pmb{u}-\dive\pmb{u}_h||_{L^{\infty}({\Omega_1\cup\Omega_2})}$, for Method 1, in Figure~\ref{fig:divergence error macro/full stabilization example 1}, and for Method 2, in Figure~\ref{fig:divergence error method2 example 1}, using $\RT_0\times Q_0$. We observe that the divergence of the velocity is deteriorated by the standard ghost penalty term for the pressure stabilization. Using a macro-element partition, we can improve the divergence error of Method 1 by having fewer elements polluted by the stabilization term $s_p$. Compare the left and the right panels of Figure~\ref{fig:divergence error macro/full stabilization example 1}.  In contrast, we see in Figure~\ref{fig:divergence error method2 example 1} that with the proposed alternative stabilization $s_b(\bfu_h,q_h)$ replacing $s_p(p_h,q_h)$, the optimal approximation of the divergence is preserved,  note that the magnitude of the error is of the order of rounding errors. 

\begin{remark}	\textbf{Preconditioner}
A simple diagonal preconditioner may be applied to great effect for the element pair $\RT_0\times Q_0$. Let $\mathbf{P}_{\pmb{u}}$ and $\mathbf{P}_p$ be the matrices associated with the following two bilinear forms,
\begin{align}\label{eq:velmassmat}
	\mathbf{P}_{\pmb{u}} &\longleftrightarrow a(\pmb{v}_h,\pmb{v}_h), \\
	\mathbf{P}_p &\longleftrightarrow (q_h,q_h)_{\Omega_1\cup\Omega_2}.
\end{align}
Then we construct the preconditioner
\begin{equation}\label{eq:precond}
	\mathbf{P} := \begin{bmatrix}
		\text{diag}(\mathbf{P}_{\pmb{u}}) & \mathbf{0} \\
			\mathbf{0} & \text{diag}(\mathbf{P}_p)
	\end{bmatrix}
\end{equation}
and consider $\mathbf{P}^{-1/2}\textbf{C}\mathbf{P}^{-1/2}\textbf{x} = \mathbf{P}^{-1/2}\textbf{b},$ where $\textbf{C}\textbf{x}=\textbf{b}$ is the linear system in equation \eqref{eq:unstabdarcy}. The solution to this preconditioned system enjoys all the desired properties and 
the condition number of the system matrix $\mathbf{P}^{-1/2}\textbf{C}\mathbf{P}^{-1/2}$ is robust with respect to how the interface cuts through the mesh, see Figure \ref{fig:example 1 precond method}. The strategy unfortunately fails for higher order elements than $\RT_0\times Q_0$, as we see in Figure \ref{fig:example 1 precond method}.  The authors are unaware of any preconditioning strategy which works in the unfitted case for higher order element pairs than $\RT_0\times Q_0$.
\begin{figure}[h]
    \centering
        \begin{subfigure}[b]{0.5 \textwidth}  	
	\centering	
    \includegraphics[scale=0.33]{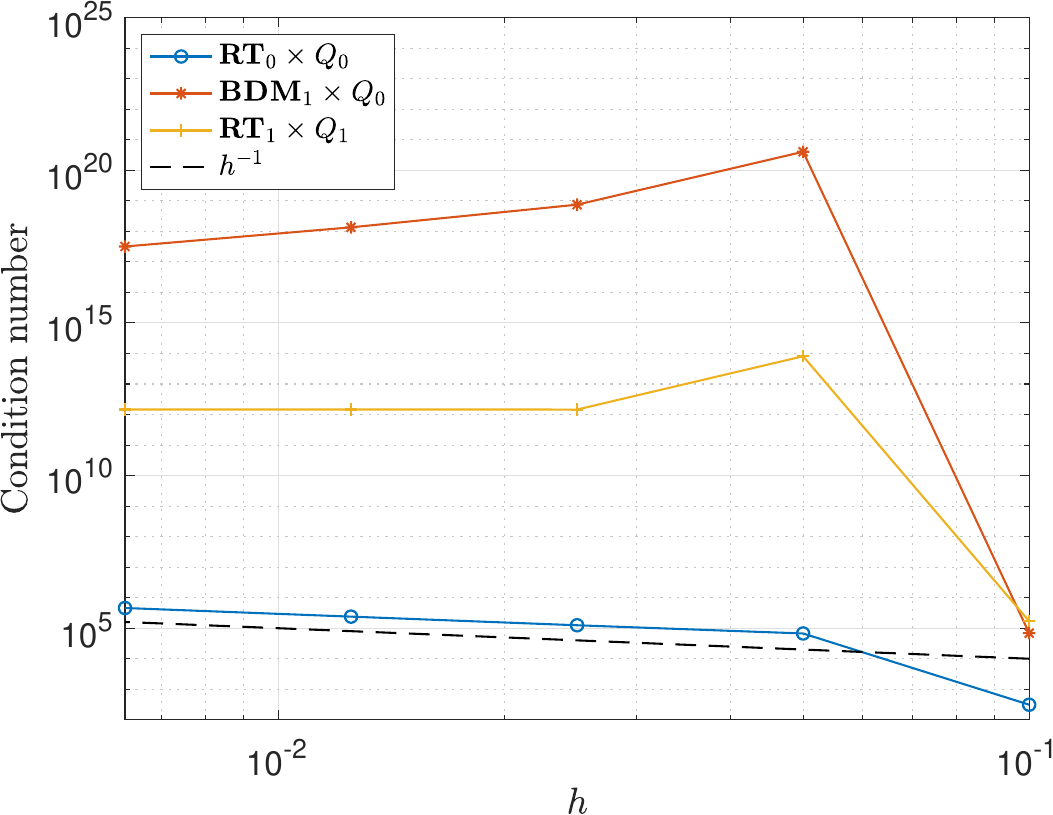}  
      \end{subfigure} 
    \caption{Example 1: The condition number of the system matrix, $\mathbf{P}^{-1/2}\textbf{C}\mathbf{P}^{-1/2}$, for the preconditioned system.}
    \label{fig:example 1 precond method}
\end{figure}
\end{remark}

\subsubsection{Extension to 3D}
We extend Example 1 to three space dimensions. The computational domain is $\Omega = [0 , 1 ]\times  [0 , 1 ] \times  [0 , 1 ]$ and the interface is the sphere defined as $$\Gamma = \left\{ (x,y,z) \in \RR^3 : r:=\sqrt{(x-1/2)^2+(y-1/2)^2+(z-1/2)^2} = R \right\},$$ with $R=0.25$. Here we have \[ g = \begin{cases} -3/R^2 \text{ in } \Omega_1,\\-6/R^2 \text{ in } \Omega_2,\end{cases}\] 
and all the other parameters are defined as in Example 1.

\begin{figure}[h!]
	\centering
	\begin{subfigure}[b]{0.45 \textwidth}  	
		\centering	
		\includegraphics[scale=0.25]{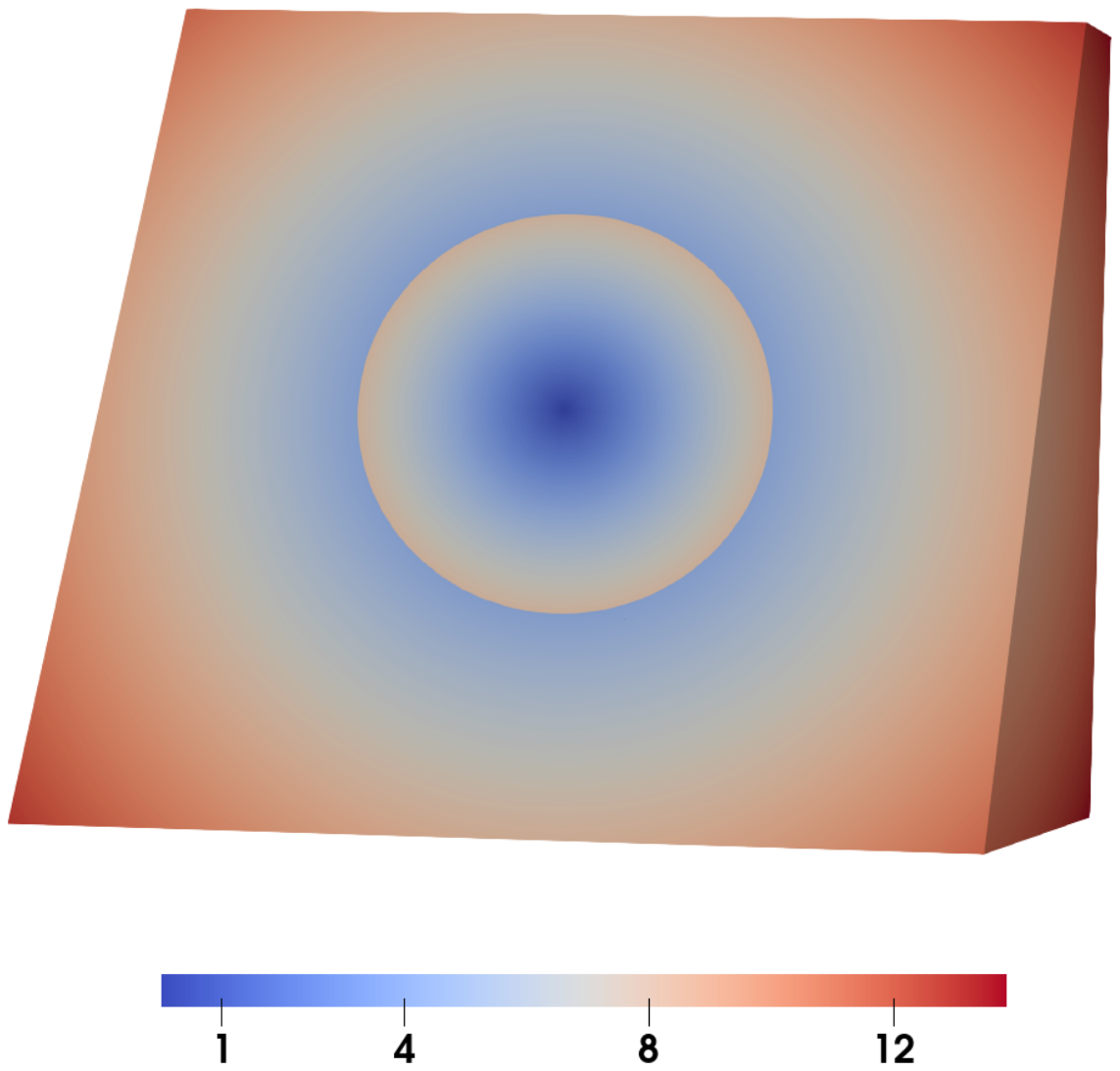}  
	\end{subfigure}
	\hfill
	\begin{subfigure}[b]{0.45\textwidth}  	
		\centering	
		\includegraphics[scale=0.25]{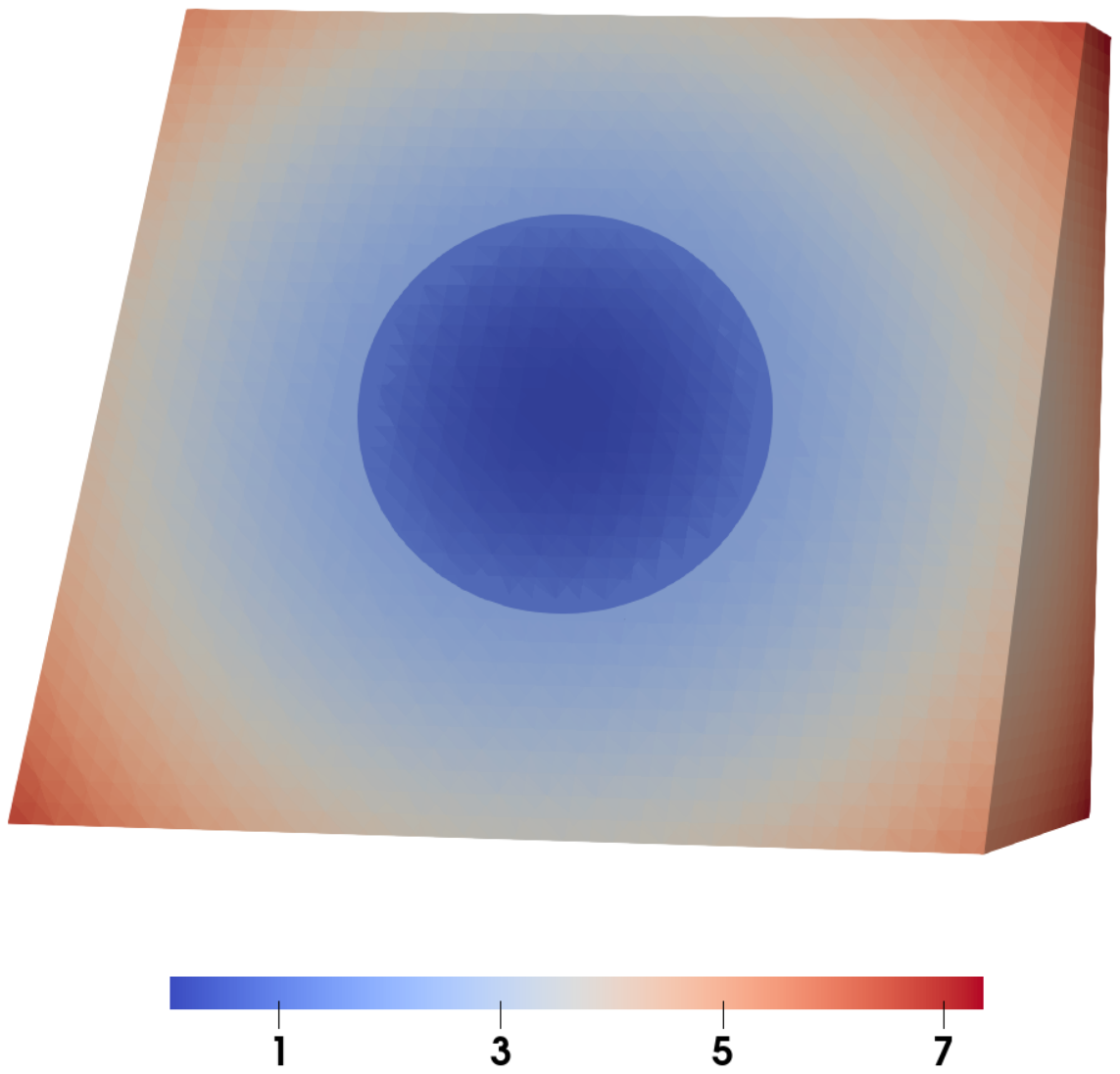}
	\end{subfigure}
	\caption{   
		Example 1 (3D): Solution obtained using Method 2 with the element pair $\RT_0\times Q_0$ and full stabilization. Left: Magnitude of the approximated velocity. Right: Approximated pressure.}
	\label{fig:example 3D solution}
\end{figure}

\begin{figure}[h!]
  	\centering
	\begin{subfigure}[b]{0.45 \textwidth}  	
		\includegraphics[scale=0.12]{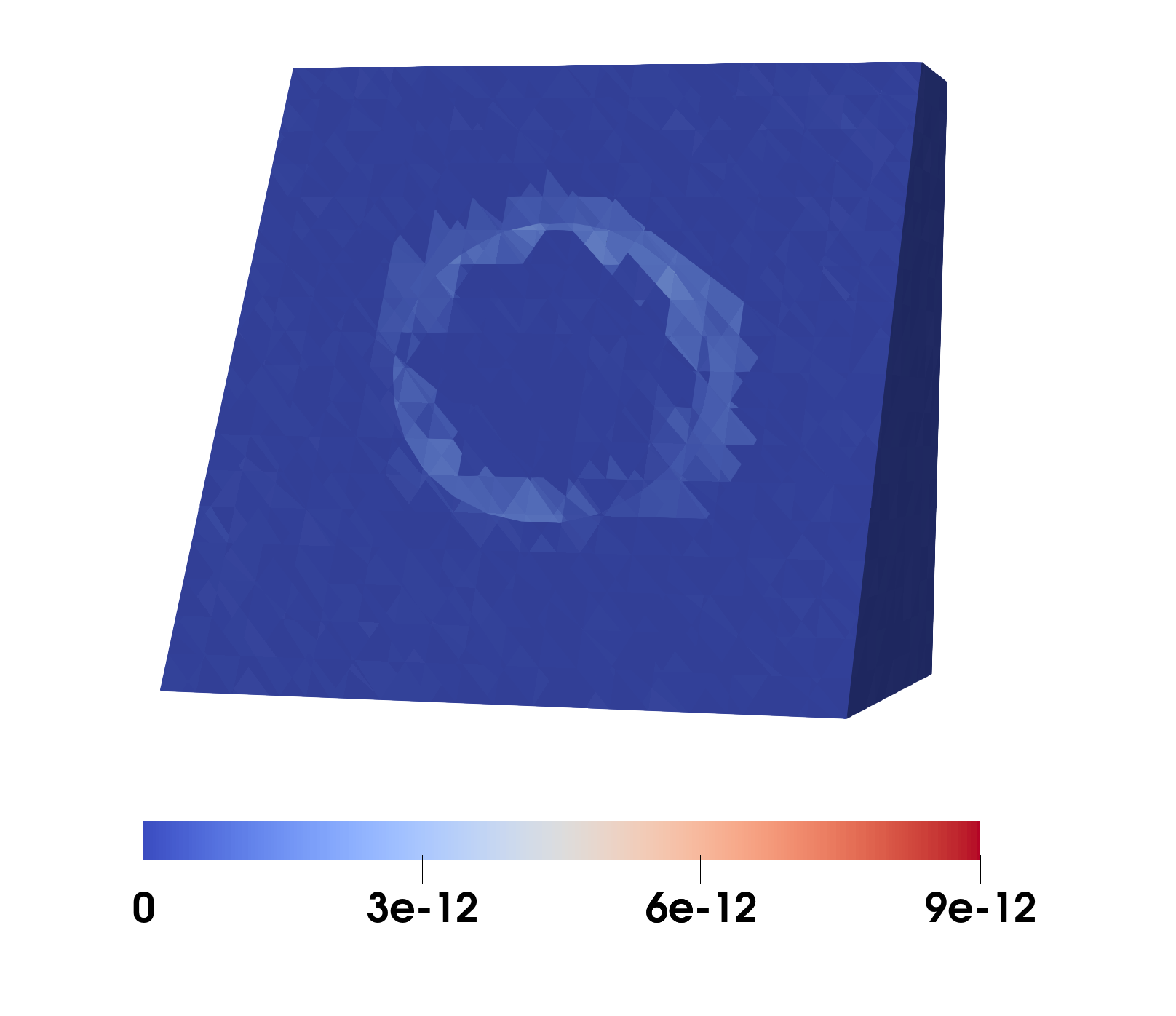}  
	\end{subfigure}
	\caption{   
		Example 1 (3D): The max-norm error $||\dive\pmb{u}-\dive\pmb{u}_h||_{L^{\infty}}$ using Method 2 with the element pair $\RT_0\times Q_0$.} 
	\label{fig:example 3D divergence error}
\end{figure}

Here we only use the element pair $\RT_0\times Q_0$. In Figure~\ref{fig:example 3D solution} we show the computed solution and in Figure~\ref{fig:example 3D divergence error} we show the divergence error in the $L^{\infty}$-norm using Method 2. Table \ref{table:3D} shows the convergence of the approximated velocity and pressure using Method 2. We obtain optimal rates of convergence and the error in the divergence is of the order of machine epsilon. 
\begin{table}
	\centering
	\begin{tabular}{|c||c||cc||cc||c||c|}
		\hline
		& h & $||p-p_h||_{L^{2}}$ & rate & $||\pmb{u}-\pmb{u}_h||_{L^{2}}$ & rate \\
		\hline
		$\RT_0\times Q_0$   & 1.00e-1 &  1.71e-1  & - 		&  2.66e-2 & -	\\
		& 5.00e-2	&  8.38e-2  & 1.02     &  8.86e-3 & 1.50	\\
		& 2.50e-2 	& 4.16e-2   & 1.01 	&  2.34e-3 & 1.89	\\
		& 1.25e-2 &  2.07e-2  & 1.00 	&  5.87e-4 & 1.99 	\\
		\hline
	\end{tabular}
	\caption{Example 1 (3D): Errors and convergence rates for pressure and velocity using Method 2 with the element pair $\RT_0\times Q_0$.}
	\label{table:3D}
\end{table}

\subsection{Example 2 ($\dive \bfu$ outside the pressure space)}
Let $g(x,y) = 100x^6+20x^2-42y^2+200y^4$ define the interface via its level set $g-c=0$ with $c=10$. Let $\bfn$ be its inward pointing normal. We choose $\xi_0=1/8$, $\eta_{\Gamma}=2c/100$, $\hat{p}=\xi_0 \eta_{\Gamma}+c/200$, $\eta=1$, and 
\begin{equation}
	p =\begin{cases}
		g/100, \ &\text{in }\Omega_1\\
		0, \ &\text{in }\Omega_2
	\end{cases},\
	\bfu = \begin{cases}
		-\nabla g / |\nabla g|, \ &\text{in }\Omega_1\\
		0, \ &\text{in }\Omega_2
	\end{cases}.
\end{equation} 

\begin{figure}[h!]
	\centering
	\includegraphics[scale=0.28]{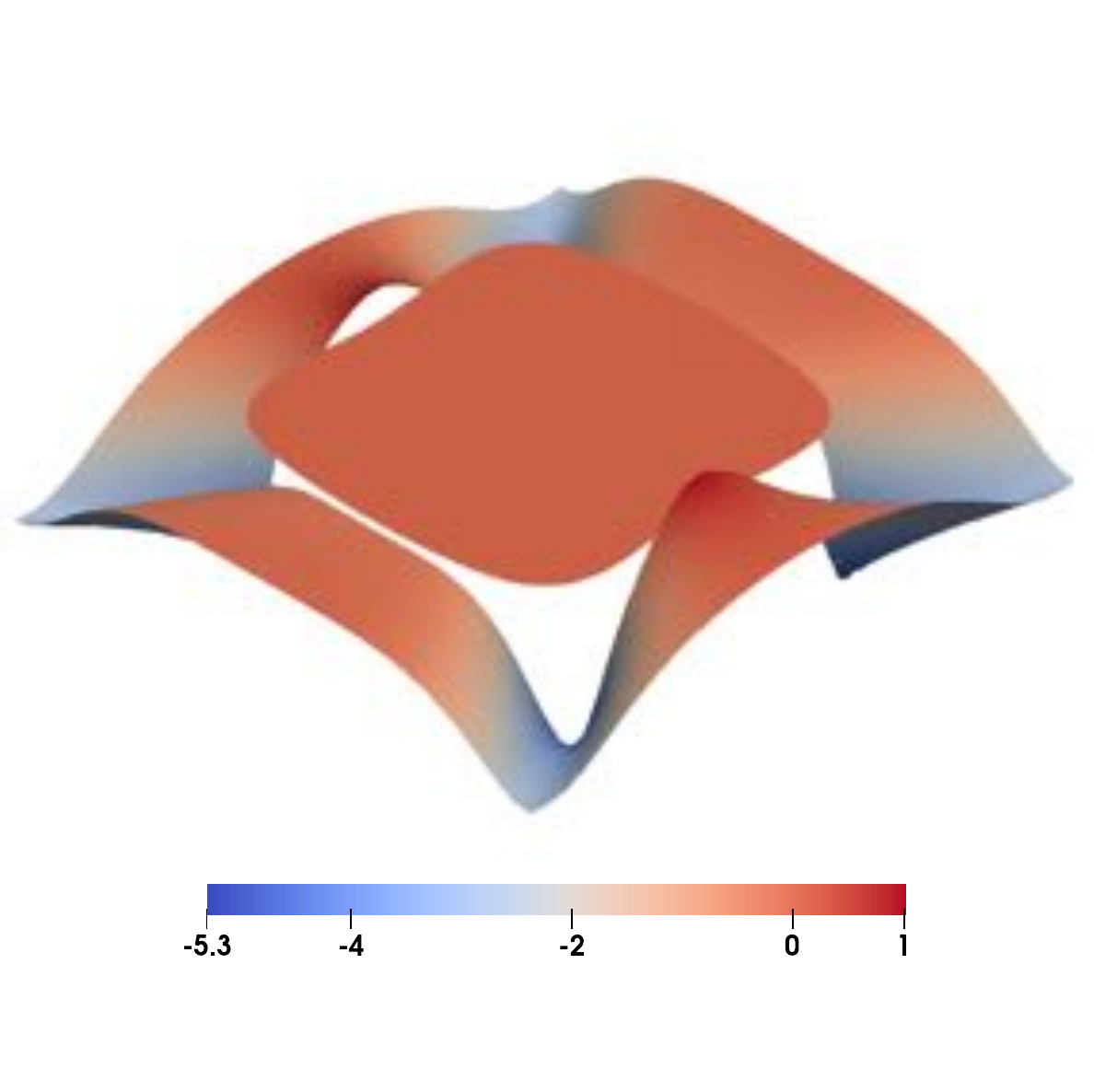}  
	\caption{Example 2: $\dive \bfu_h$ obtained using Method 2 with the element pair $\RT_1\times Q_1$ and a macro-element partitioning with $\delta = 0.25$.}
	\label{fig:example 2 divergence}
\end{figure}

The interface looks like a saddle and has both positive and negative curvature, see Figure~\ref{fig:example 2 divergence} where the divergence $\dive\bfu$ is shown in the domain $\Om$. Of particular note is that the exact divergence $\dive \bfu$ is not in the pressure space for any of the considered element pairs. 
Thus, the best we can hope for in terms of the divergence error is to match the behaviour of standard FEM where the mesh is fitted to the interface. Here we consider the pure Dirichlet case of $\boundp=\emptyset$, and $\bfu|_{\partial\Omega}=-\nabla g/|\nabla g|.$

\begin{figure}[h]
\centering 
\includegraphics[scale=0.53]{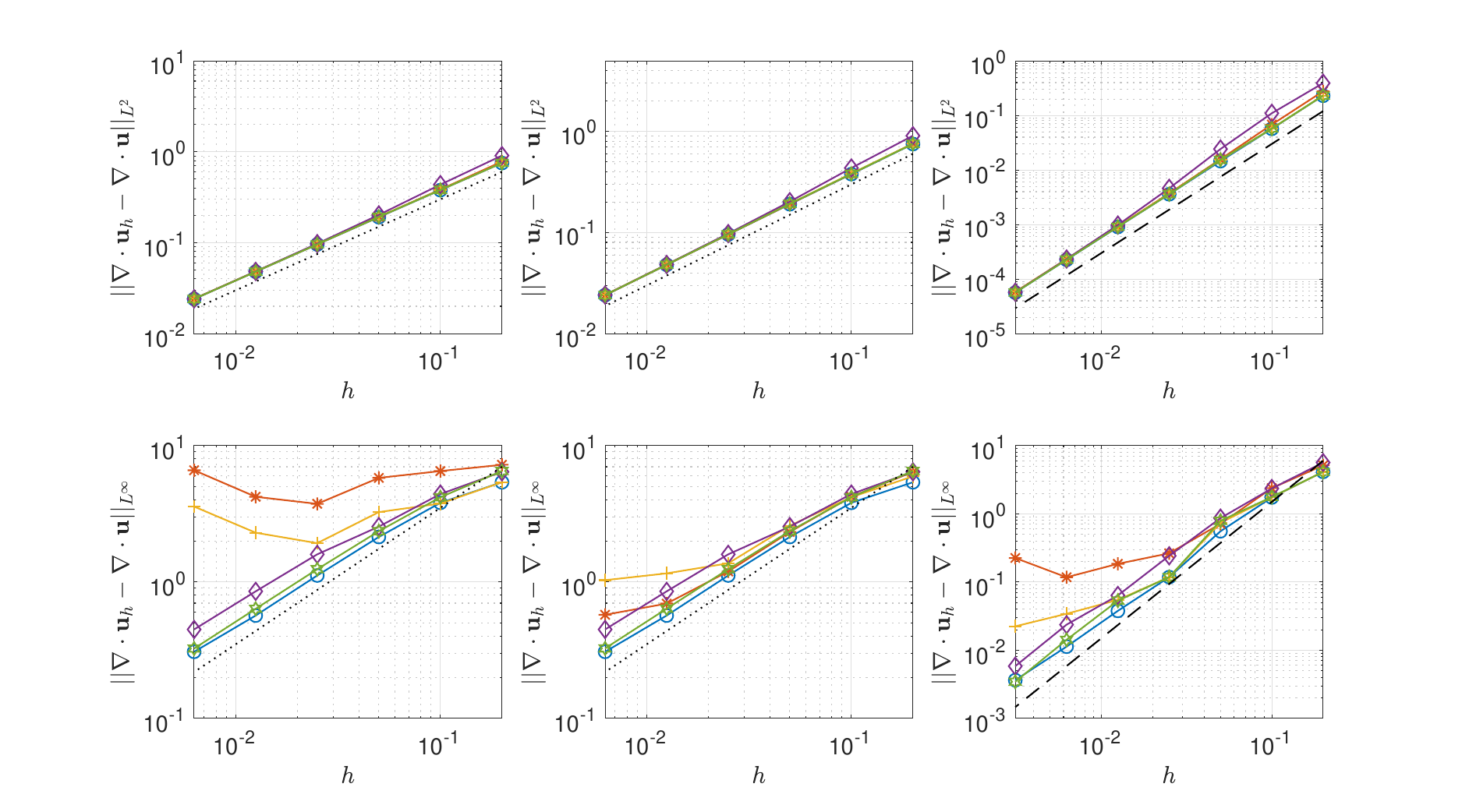} \\
\includegraphics[scale=0.35]{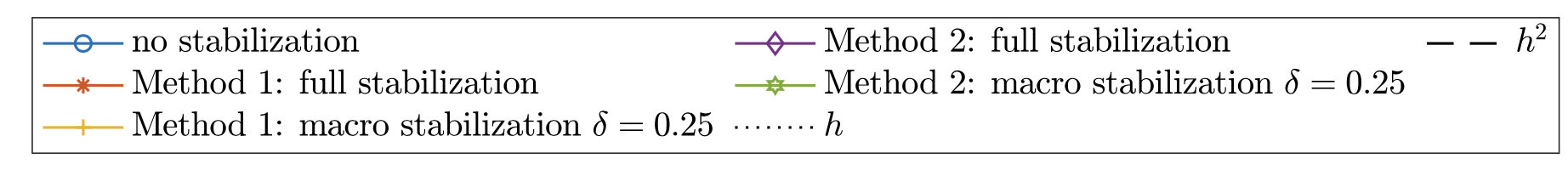}
\caption{Example 2:  The L$^2$-error (first row) and pointwise error (second row) of the divergence versus mesh size h. Left: $\RT_0\times Q_0$. Middle:$\BDM_1\times Q_0$. Right: $\RT_1\times Q_1$.}
 \label{fig:plot example2 convergence divergence}  	
\end{figure}

In Figure~\ref{fig:plot example2 convergence divergence} we show how the presented methods differ in the way that the divergence of the computed solution converge. Method 2 and the unstabilized method match standard FEM and converge optimally in the $L^2-$ and $L^{\infty}-$norms, while Method 1 does not converge at all in the $L^{\infty}$-norm. We also see that using a macro-element partitioning and stabilizing less reduces the constant in the $L{^\infty}$-error. Method 2 with macro stabilization has for all element pairs errors of similiar magnitude as the unstabilized method. 
\begin{figure}[h]
\centering 
\hspace*{-2cm}   
\includegraphics[scale=0.58]{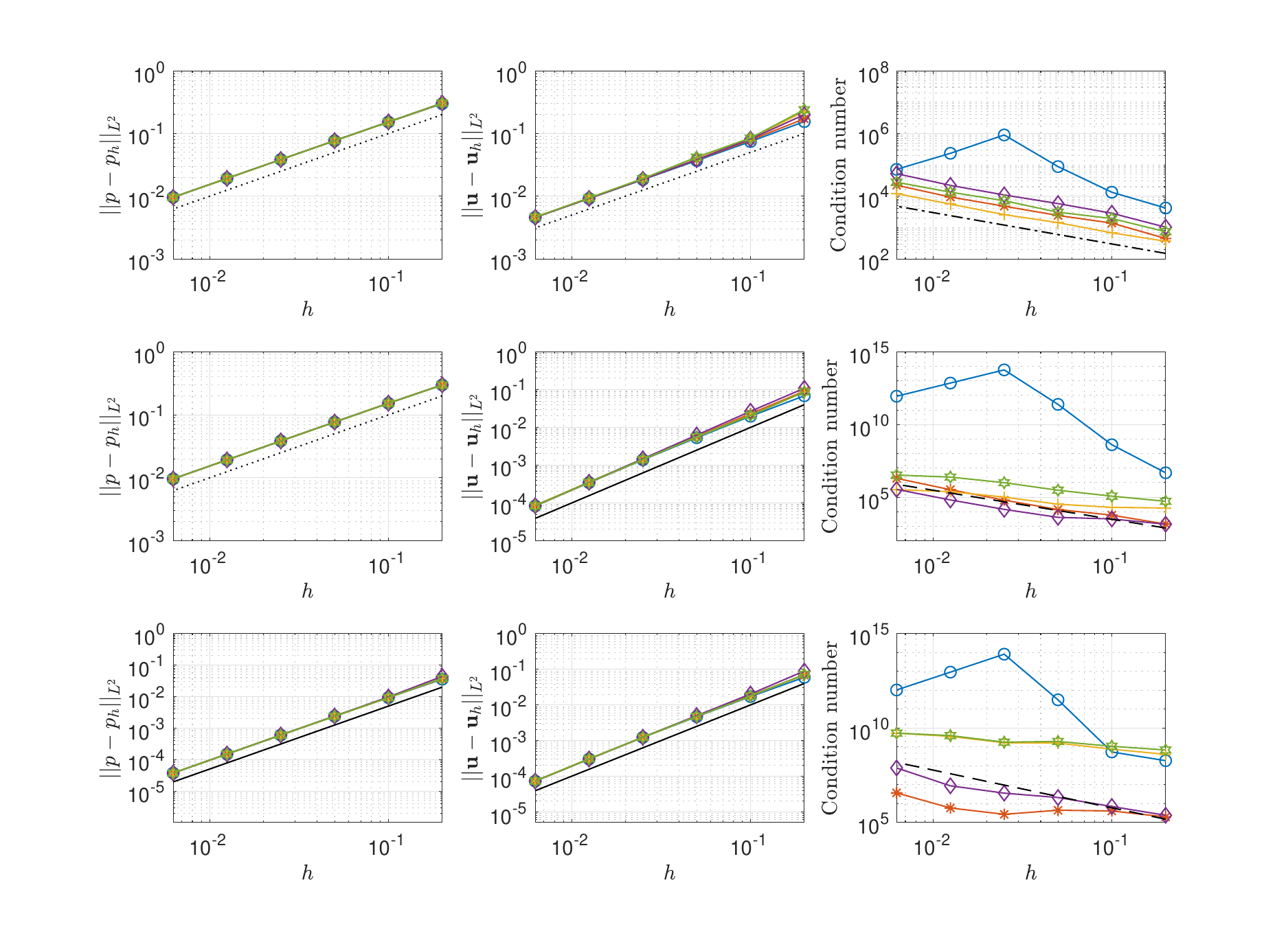} \\
\vspace{-0.5cm}
\includegraphics[scale=0.35]{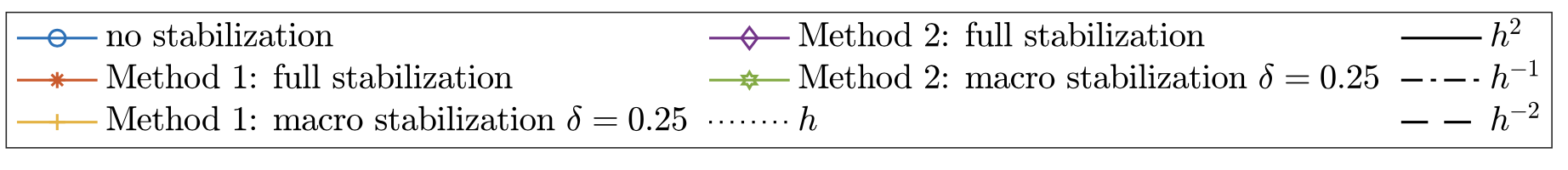}
\caption{Example 2: Results using element pairs $\RT_0\times Q_0$ (first row), $\BDM_1\times Q_0$ (second row) and $\RT_1\times Q_1$ (third row). 
        Left: The $L^2$-error of the pressure versus mesh size h. Middle:  The $L^2$-error of the velocity field versus mesh size h. Right: The spectral condition number versus mesh size h.}     
 \label{fig:plot example2 convergence}  	
\end{figure}
In Figure~\ref{fig:plot example2 convergence} we show the $L{^2}$-error of the approximated velocity and pressure as well as the spectral condition number as function of mesh size. We observe optimal convergence rates; first order for both velocity and pressure using $\RT_0\times Q_0$, second order for the velocity and first order for the pressure using $\BDM_1\times Q_0$, and second order convergence both for the velocity and pressure using $\RT_1\times Q_1$.  
        
\section{Conclusion}\label{sec:conclusion}
We studied two cut finite element discretizations for Darcy flow in fractured porous media, where the fracture is modeled by an interface (an internal boundary). Both methods use stabilization terms in the variational formulation in order to provide well-posed discretizations with accurate velocity and pressure approximations independently of how the interface cuts through the computational mesh. However, the first cut finite element method, Method 1, which is  based on standard ghost penalty terms,  one for velocity $s_{\bfu}(\bfu_h,\bfv_h)$ and one for pressure $s_p(p_h,q_h)$, fails in preserving the good properties of Raviart-Thomas and Brezzi-Douglas-Marini element pairs. Therefore, we developed an alternative stabilization that preserves the divergence-free property of these elements also in the unfitted setting. With the new stabilization term the cut finite element discretization we propose, Method 2, keeps the saddle point structure and provides a good alternative to the classical fitted $\bfH^{\dive}$ conforming set up. We also used so called macro-element stabilization, where stabilization is applied very restrictively,  instead of full stabilization. This resulted in that pointwise errors in the divergence, using the proposed method, were similar to errors from the unstabilized discretization. We believe this work is an important step in the development of divergence-free approximations in the unfitted setting. The new method and the analysis we present in this paper hold for the case when the mesh is fitted to the outer boundary and boundary conditions for the velocity are enforced strongly. The interface is however arbitrarily positioned relative to the mesh. Our aim is to consider the case of an unfitted outer boundary in a forthcoming paper.
                              
\section*{Acknowledgments}
The authors are grateful to the anonymous reviewers for carefully reading our manuscript and providing valuable comments that helped us to improve this work.
\addcontentsline{toc}{section}{Bibliography}
\bibliographystyle{siamplain}                           
\bibliography{ref} 
     
\end{document}